\documentclass[12pt]{amsart}
\usepackage{latexsym, amssymb, amsmath, amscd}
 \usepackage{eucal}

\theoremstyle{plain}
    \newtheorem{rema}{Remark}[section]
    \newtheorem{propo}[rema]{Proposition}

   \newtheorem{theo}[rema]{Theorem}
 \newtheorem{conj}[rema]{Conjecture}
   \newtheorem{defi}[rema]{Definition}
    \newtheorem{lemma}[rema]{Lemma}
    \newtheorem{corol}[rema]{Corollary}
     \newtheorem{exam}[rema]{Example}
  \newtheorem{rmk}[rema]{Remark}

	\newcommand{\nno}{\nonumber}

	\newcommand{\p}{\partial}
\newcommand{\fr}{\frac}

 \newcommand{\pf}{{\it Proof:}\hspace{2ex}}
 
 \newcommand{\epfv}{\hspace{1em}$\Box$\vspace{1em}}

\newcommand{\bR}{{\mathbb R}}

\newcommand{\bN}{{\mathbb N}}

\newcommand{\BQ}{\begin{eqnarray}}
\newcommand{\EQ}{\end{eqnarray}}
\newcommand{\BQn}{\begin{eqnarray*}}
\newcommand{\EQn}{\end{eqnarray*}}

\newcommand{\BL}{\begin{align}}
\newcommand{\EL}{\end{align}}
\newcommand{\BLn}{\begin{align*}}
\newcommand{\ELn}{\end{align*}}

\newcommand{\Hes}{ \text{Hes\,} }

\newcommand{\bC}{\mathbb C}

\newcommand{\bCz}{\bC[z]}

\newcommand{\cA}{{\mathcal A}}

\newcommand{\bD}{{\mathbb D}}
\newcommand{\cD}{{\mathcal D}}

\newcommand{\JC}{\mbox{\bf JC}}
\newcommand{\VC}{\mbox{\bf VC}}
\newcommand{\HVC}{\mbox{\bf HVC}}
\newcommand{\la}{\langle}
\newcommand{\ra}{\rangle}

\renewcommand{\theequation}{\thesection.\arabic{equation}}
\renewcommand{\therema}{\thesection.\arabic{rema}}
\newcommand{\EAn}{\end{align*}}
\newcommand{\wtilde}{\widetilde}

\title[A Vanishing Conjecture on Differential Operators]
{A Vanishing Conjecture on Differential Operators with Constant 
Coefficients}

    \author{Wenhua Zhao}      
   \date{\today}
    
\begin{document}

\begin{abstract}
In the recent progress \cite{BE1}, \cite{Me} and \cite{HNP}, 
the well-known JC (Jacobian conjecture) (\cite{BCW}, \cite{E}) 
has been reduced to a VC (vanishing conjecture) on 
the Laplace operators and HN (Hessian nilpotent) polynomials 
(the polynomials whose Hessian matrix are nilpotent).
In this paper, we first show that the vanishing conjecture above, 
hence also the JC,  is equivalent to a vanishing conjecture 
for all 2nd order homogeneous differential operators $\Lambda$ 
and $\Lambda$-nilpotent polynomials $P$ 
(the polynomials $P(z)$ satisfying 
$\Lambda^m P^m=0$ for all $m\ge 1$).
We then transform some results in the literature
on the JC, HN polynomials and the VC 
of the Laplace operators to certain results on 
$\Lambda$-nilpotent polynomials 
and the associated VC for 2nd order 
homogeneous differential operators $\Lambda$.
This part of the paper can also be read 
as a short survey on HN polynomials 
and the associated VC in the more general setting. 
Finally, we discuss a still-to-be-understood
connection of $\Lambda$-nilpotent 
polynomials in general with the classical 
orthogonal polynomials in one or more variables. 
This connection provides a conceptual 
understanding for the isotropic properties 
of homogeneous $\Lambda$-nilpotent 
polynomials for 2nd order homogeneous full rank
differential operators $\Lambda$ with 
constant coefficients.
\end{abstract}

   

 \bibliographystyle{alpha}
    \maketitle


\renewcommand{\theequation}{\thesection.\arabic{equation}}
\renewcommand{\therema}{\thesection.\arabic{rema}}
\setcounter{equation}{0}
\setcounter{rema}{0}
\setcounter{section}{0}

\section{\bf Introduction}\label{S1}

Let $z=(z_1, z_2, \dots, z_n, \dots)$ 
be a sequence of free commutative variables 
and $D=(D_1, D_2, \dots, D_n, \dots)$ with 
$D_i\!:=\frac{\p}{\p z_i}$ $(i\geq 1)$. 
For any $n\geq 1$, denote by $\cA_n$ (resp.\,$\bar{\cA}_n$) 
the algebra of polynomials  
(resp.\,formal power series) in $z_i$ 
$(1\leq i\leq n)$. Furthermore, 
we denote by $\cD[\cA_n]$ or $\cD[n]$
(resp.\,$\bD[\cA_n]$ or $\bD[n]$) the algebra of 
differential operators of the polynomial 
algebra $\cA_n$ (resp.\,with constant coefficients). 
Note that, for any $k\geq n$, elements of $\cD[n]$ 
are also differential operators of $\cA_k$ and  
$\bar{\cA}_k$. For any $d\geq 0$, denote by 
$\bD_d[n]$ the set of {\it homogeneous 
differential operators of order 
$d$ with constants coefficients}. 
We let $\cA$ (resp.\,${\bar\cA}$) 
be the union of $\cA_n$ (resp.\,$\bar {\cA}_n$) 
$(n\geq 1)$, $\cD$ (resp.\,$\bD$) 
the union of $\cD[n]$ (resp.\,$\bD[n]$) $(n\geq 1)$, 
and, for any $d\geq 1$, $\cD_d$ the 
union of $\cD_d[n]$ $(n\geq 1)$.

Recall that $\JC$ (the Jacobian conjecture) 
which was first proposed 
by Keller \cite{Ke} in 1939, 
claims that, {\it for any polynomial map $F$ 
of $\bC^n$ with Jacobian $j(F)=1$, 
its formal inverse map $G$ must also be a polynomial map}. 
Despite intense study 
from mathematicians in more than sixty years, 
the conjecture is still open
even for the case $n=2$. 
For more history and known results 
before $2000$ on $\JC$, 
see \cite{BCW}, \cite{E} 
and references there.

Based on the remarkable symmetric reduction  achieved 
in \cite{BE1}, \cite{Me} and the classical celebrated 
homogeneous reduction \cite{BCW} and \cite{Y} 
on $\JC$, the author in \cite{HNP}
reduced $\JC$  further to the following 
{\it vanishing conjecture} on the Laplace operators 
$\Delta_n\!:=\sum_{i=1}^n D_i^2$ of the polynomial algebra 
$\cA_n$ and {\it HN} ({\it Hessian nilpotent}) 
polynomials $P(z)\in \cA_n$, 
where we say a polynomial or formal power series 
$P(z)\in \bar{\cA}_n$ is HN if its 
Hessian matrix 
$\Hes (P)\!:=(\frac{\p^2 P}{\p z_i\p z_j})_{n\times n}$ 
is nilpotent. 

\begin{conj} \label{LVC} 
For any HN $($homogeneous$)$ polynomial $P(z)\in \cA_n$ 
$($of degree $d=4$$)$, 
we have $\Delta_n^m P^{m+1}(z)=0$ 
when $m>>0$.
\end{conj}

Note that, the following criteria of Hessian nilpotency 
were also proved in Theorem $4.3$, \cite{HNP}. 

\begin{theo}\label{Crit-1}
For any $P(z)\in \bar {\cA}_n$ with $o(P(z))\geq 2$, 
the following statements are equivalent.
\begin{enumerate}
\item[(1)] $P(z)$ is HN.
\item[(2)] $\Delta^m P^m=0$ for any $m\geq 1$.
\item[(3)] $\Delta^m P^m=0$ for any $1\leq m\leq n$.
\end{enumerate}
\end{theo}

Through the criteria in the proposition above, 
Conjecture \ref{LVC} can be generalized to
other differential operators as follows 
(see Conjecture \ref{GVC} below). 

First let us fix the following notion 
that will be used throughout the paper.

\begin{defi}\label{Def-LN}
Let $\Lambda \in \cD[\cA_n]$ and $P(z) \in \bar{\cA}_n$. 
We say $P(z)$ is {\it $\Lambda$-nilpotent} 
if $\Lambda^m P^m =0$ for any $m\geq 1$. 
\end{defi}

Note that,\label{rmk-HN} when $\Lambda$ is the Laplace operator 
$\Delta_n$, by Theorem \ref{Crit-1}, 
a polynomial or formal power series $P(z)\in \cA_n$ 
is $\Lambda$-nilpotent iff it is HN. 

With the notion above, Conjecture \ref{LVC} has the 
following natural generalization to differential 
operators with constant coefficients.

\begin{conj}\label{GVC}
For any $n\geq 1$ and $\Lambda \in \bD[n]$, if $P(z)\in \cA_n$ 
is $\Lambda$-nilpotent, then $\Lambda^m P^{m+1}=0$ 
when $m>>0$.
\end{conj}

We call the conjecture above the 
{\it vanishing conjecture} for differential 
operators with constant coefficients
and denote it by $\VC$. The special case 
of $\VC$ with $P(z)$ homogeneous is called the 
{\it homogeneous vanishing conjecture} 
and denoted by $\HVC$.
When the number $n$ of variables 
is fixed, $\VC$ (resp.\,$\HVC$) is called
(resp.\,{\it homogeneous}) 
{\it vanishing conjecture in $n$ variables} 
and denoted by $\VC[n]$ (resp.\,$\HVC[n]$).

Two remarks on $\VC$ are as follows.
First, due to a counter-example given by M. de Bondt 
(see example \ref{Bondt}), $\VC$ does 
not hold in general for differential operators 
with non-constant coefficients. 
Secondly, one may also allow $P(z)$ in 
$\VC$ to be any $\Lambda$-nilpotent 
formal power series. No counter-example 
to this more general $\VC$ is known yet. 

In this paper, we first apply certain 
linear automorphisms and Lefschetz's principle
to show Conjecture \ref{LVC}, hence also $\JC$, 
is equivalent to $\VC$ or $\HVC$ for all 2nd order 
homogeneous differential operators $\Lambda\in \bD_2$ 
(see Theorem \ref{GVC-MainThm}). 
We then in Section \ref{S3} transform some results on $\JC$, 
HN polynomials and Conjecture \ref{LVC}  
obtained in \cite{Wa}, \cite{BE2}, \cite{BE3}, 
\cite{HNP}, \cite{OP-HNP} and \cite{EZ} to certain 
results on $\Lambda$-nilpotent 
$(\Lambda\in \bD_2)$ polynomials 
and $\VC$ for $\Lambda$. Another purpose 
of this section is to give a survey on recent study 
on Conjecture \ref{LVC} and HN polynomials in the more
general setting of $\Lambda\in \bD_2$ and 
$\Lambda$-nilpotent polynomials. 
This is also why some results in the general setting,
even though their proofs are straightforward, are 
also included here.

Even though, due to M. de Bondt's counter-example 
(see Example \ref{Bondt}), $\VC$ does not hold for 
all differential operators with non-constant coefficients, 
it is still interesting to consider whether or not $\VC$ holds
for higher order differential operators 
with constant coefficients; and if it also holds   
even for certain families of differential operators 
with non-constant coefficients. For example, 
when $\Lambda=D^{\bf a}$ with ${\bf a}\in \bN^n$ 
and $|{\bf a}|\ge2$, $\VC[n]$ for $\Lambda$ is equivalent 
to a conjecture on Laurent polynomials 
(see Conjecture \ref{LP-Conj}). This conjecture 
is very similar to a non-trivial theorem 
(see Theorem \ref{DK1}) on Laurent polynomials, 
which was first conjectured  by O. Mathieu \cite{Mat} 
and later proved by 
J. Duistermaat and W. van der Kallen \cite{DK}.

In general, to consider the questions above, 
one certainly needs to 
get better understandings on the 
$\Lambda$-nilpotency condition, 
i.e. $\Lambda^m P^m=0$ for any $m\geq 1$.  
One natural way to look at this condition 
is to consider the sequences 
of the form $\{ \Lambda^m P^m\,|\, m\geq 1\}$
for general differential operators $\Lambda$ and 
polynomials $P(z)\in \cA$. What special 
properties do these sequences have 
so that $\VC$ wants them all vanish? 
Do they play any important roles in 
other areas of mathematics? 

The answer to the first question above 
is still not clear. 
The answer to the later seems to be ''No". 
It seems that the sequences of the form 
$\{ \Lambda^m P^m\,|\, m\geq 1\}$
do not appear very often in mathematics. 
But the answer turns out to be ``Yes" 
if one considers the question in the setting 
of some localizations $\mathcal B$ of $\cA_n$. 
Actually, as we will discuss in some detail in 
subsection  \ref{S4.1}, all classical 
orthogonal polynomials \label{nil-orp1}
in one variable have the form  
$\{ \Lambda^m P^m\,|\, m\geq 1\}$
except there one often chooses 
$P(z)$ from some localizations 
$\mathcal B$ of $\cA_n$ and $\Lambda$ 
a differential 
operators of $\mathcal B$. Some classical 
polynomials in several variables can also 
be obtained from sequences of the form 
$\{ \Lambda^m P^m\,|\, m\geq 1\}$
by a slightly modified procedure. 

Note that, due to their applications in many different 
areas of mathematics, especially in ODE, PDE, the
eigenfunction problems and representation theory, 
orthogonal polynomials have been under intense 
study by mathematicians in the last two centuries. 
For example, in \cite{SHW} 
published in $1940$, about $2000$ 
published articles mostly 
on one-variable orthogonal polynomials  
have been included. The classical reference 
for one-variable orthogonal polynomials 
is \cite{Sz} (see also \cite{AS}, 
\cite{C}, \cite{Si1}). For multi-variable 
orthogonal polynomials, 
see \cite{DX}, \cite{Ko} 
and references there.

It is hard to believe that the connection 
discussed above between $\Lambda$-nilpotent polynomials 
or formal power series and classical orthogonal polynomials 
is just a coincidence. But a precise understanding 
of this connection still remains mysterious. 
What is clear is that, $\Lambda$-nilpotent polynomials 
or formal power series and the polynomials 
or formal power series $P(z)\in \bar{\cA}_n$ 
such that the sequence $\{ \Lambda^m P^m\,|\, m\geq 1\}$
for some differential operator $\Lambda$ provides 
a sequence of orthogonal polynomials 
lie in two opposite extreme sides, since, 
from the same sequence 
$\{ \Lambda^m P^m\,|\, m\geq 1\}$,
the former provides nothing but zero;
while the later provides an orthogonal  
basis for $\cA_n$.

Therefore, one naturally expects that
$\Lambda$-nilpotent polynomials 
$P(z)$$\in \cA_n$ should be isotropic \label{nil-orp2}
with respect to a certain $\bC$-bilinear form 
of $\cA_n$. It turns out that, as we will show 
in Theorem \ref{A-Isotropic} and Corollary \ref{A-Isotropic-C}, it 
is indeed the case when $P(z)$ is homogeneous and 
$\Lambda\in \bD_2[n]$ is of full rank. Actually,
in this case $\Lambda^m P^{m+1}$ $(m\ge 0)$  
are all isotropic with respect to same 
properly defined $\bC$-bilinear form. 
Note that, Theorem \ref{A-Isotropic} 
and Corollary \ref{A-Isotropic-C}
are just transformations of the isotropic properties of 
HN nilpotent polynomials, which were first 
proved in \cite{HNP}. But the proof in \cite{HNP} 
is very technical and lacks any convincing 
interpretations. 
From the ``formal" connection of 
$\Lambda$-nilpotent polynomials 
and orthogonal polynomials discussed above, 
the isotropic properties of homogeneous 
$\Lambda$-nilpotent polynomials with 
$\Lambda\in\bD_2[n]$ of full rank become 
much more natural.

The arrangement of the paper is as follows.
In Section \ref{S2}, we mainly show that 
Conjecture \ref{LVC}, hence also $\JC$, 
is equivalent to $\VC$ or $\HVC$ 
for all $\Lambda \in \bD_2$ 
(see Theorem \ref{GVC-MainThm}).
One consequence of this equivalence
is that, to prove or disprove $\VC$ or $\JC$, 
the Laplace operators are not the 
only choices, even though they 
are the best in many situations.
Instead, one can choose any sequence 
$\Lambda_{n_k} \in \bD_2$
with strictly increasing ranks 
(see Proposition \ref{MainPropo-3}).
For example, one can choose 
the ``Laplace operators" with respect to 
the Minkowski metric or symplectic metric, 
or simply choose $\Lambda$ to be
the complex $\bar{\p}$-Laplace operator
$\Delta_{\bar{\p}, k}$ $(k\ge 1)$
in Eq.\,$(\ref{cx-Delta})$.

In Section \ref{S3}, we transform some results 
on $\JC$, HN polynomials and Conjecture \ref{LVC} 
in the literature to certain results 
on $\Lambda$-nilpotent $(\Lambda\in \bD_2)$ 
polynomials $P(z)$ and $\VC$ for $\Lambda$. 
In subsection  \ref{S3.1}, 
we discuss some results on the polynomial maps 
and PDEs associated with $\Lambda$-nilpotent 
polynomials for $\Lambda\in \bD_2[n]$ of full rank 
(see Theorems \ref{PDE1}--\ref{Heat}). The results 
in this subsection are transformations of 
those in \cite{Burgers} and \cite{HNP} 
on HN polynomials and their 
associated symmetric polynomial 
maps.

In subsection  \ref{S3.2}, we give four criteria 
of $\Lambda$-nilpotency  $(\Lambda\in \bD_2)$ 
(see Propositions \ref{A-Crit-0}, \ref{A-Crit-1}, 
\ref{A-Crit-3} and \ref{A-Crit-4}). The criteria 
in this subsection are transformations of the 
criteria of Hessian nilpotency 
derived in \cite{HNP} and \cite{OP-HNP}.
In subsection  \ref{S3.3}, we transform
some results in \cite{BCW}, \cite{Wa} and \cite{Y} 
on $\JC$; \cite{BE2} and \cite{BE3} on 
symmetric polynomial maps; \cite{HNP}, \cite{OP-HNP} 
and \cite{EZ} on HN polynomials to certain results on $\VC$ 
for $\Lambda\in \bD_2$. Finally, 
we recall a result in \cite{OP-HNP} which says,    
$\VC$ over fields $k$ of characteristic $p>0$,
even under some conditions weaker than 
$\Lambda$-nilpotency, actually holds 
for any differential operators $\Lambda$ 
of $k[z]$ (see Proposition \ref{pVC} 
and Corollary \ref{pVC-C1}).

In subsection  \ref{S3.4}, we consider $\VC$ 
for high order differential operators with 
constant coefficients. In particular, 
we show in Proposition \ref{bD-1k} $\VC$ holds 
for $\Lambda=\delta^k$ $(k\geq 1)$, 
where $\delta$ is a derivation 
of $\cA$. In particular, $\VC$ holds 
for any $\Lambda\in \bD_1$ 
(see Corollary \ref{bD1}). We also show that, 
when $\Lambda=D^{\bf a}$ with ${\bf a}\in \bN^n$ 
and $|{\bf a}|\ge 2$, $\VC$ is equivalent to
a conjecture, Conjecture \ref{LP-Conj}, 
on Laurent polynomials. This conjecture 
is very similar to a non-trivial 
theorem (see Theorem \ref{DK1}) 
first conjectured by O. Mathieu  \cite{Mat} 
and later proved by J. Duistermaat 
and W. van der Kallen \cite{DK}.

In subsection  \ref{S4.1}, by using 
Rodrigues' formulas Eq.\,$(\ref{Rodrigues})$, 
we show that all classical 
orthogonal polynomials  in one variable 
have the form $\{ \Lambda^m P^m\,|\, m\geq 1\}$
for some $P(z)$ in certain localizations 
$\mathcal B$ of $\cA_n$ and $\Lambda$ 
a differential operators of $\mathcal B$.
We also show that some classical 
polynomials in several variables can also 
be obtained from sequences of the form 
$\{ \Lambda^m P^m\,|\, m\geq 1\}$ with 
a slight modification.  
Some of the most classical orthogonal polynomials 
in one or more variables are briefly 
discussed in Examples \ref{HP}--\ref{GP}, 
\ref{OPB} and \ref{OPSX}. 
In subsection  \ref{S4.2}, we transform the 
isotropic properties of homogeneous 
HN homogeneous polynomials derived in \cite{HNP} 
to homogeneous $\Lambda$-nilpotent 
polynomials for $\Lambda\in\bD_2[n]$ 
of full rank (see Theorem \ref{A-Isotropic} and 
Corollary \ref{A-Isotropic-C}).

{\bf Acknowledgment:} The author is very grateful to Michiel de Bondt 
for sharing his counterexample (see Example \ref{Bondt}) with the author,
and to Arno van den Essen for inspiring personal communications. The author 
would also like to thank the referee very much for many valuable suggestions 
to improve the first version of the paper.

\renewcommand{\theequation}{\thesection.\arabic{equation}}
\renewcommand{\therema}{\thesection.\arabic{rema}}
\setcounter{equation}{0}
\setcounter{rema}{0}

\section{\bf The Vanishing Conjecture for the 2nd Order Homogeneous
Differential Operators with Constant Coefficients}\label{S2}

In this section, we apply certain 
linear automorphisms and Lefschetz's principle
to show Conjecture \ref{LVC}, hence also $\JC$, 
is equivalent to $\VC$ or $\HVC$ for all 
$\Lambda\in \bD_2$ (see Theorem \ref{GVC-MainThm}). 
In subsection  \ref{S2.1}, we fix some notation 
and recall some lemmas that will be needed throughout this paper.
In subsection  \ref{S2.2}, we prove the main results 
of this section, Theorem \ref{GVC-MainThm} and 
Proposition \ref{MainPropo-3}.

\subsection{Notation and Preliminaries}\label{S2.1}

Throughout this paper, unless stated otherwise, 
we will keep using the notations and terminology 
introduced in the previous section and 
also the ones fixed as below.  \\

\begin{enumerate}

\item For any $P(z)\in \cA_n$, we denote by 
$\nabla P$ the {\it gradient} of $P(z)$, i.e. we set
\begin{align}
\nabla P(z)\!:=(D_1 P, \, D_2 P, \dots, D_n P).
\end{align}

\item For any $n\ge 1$, we let 
$SM(n, \bC)$ (resp.\,$SGL(n, \bC)$) denote 
the symmetric complex $n\times n$ 
(resp.\,invertible) matrices.

\item For any $A=(a_{ij}) \in SM(n, \bC)$, we set 
\begin{align}
\Delta_A\!:=\sum_{i,j=1}^n a_{ij} D_i D_j \in \bD_2[n].
\end{align}

Note that, for any $\Lambda \in \bD_2[n]$, there exists 
a unique $A \in SM(n, \bC)$ such that 
$\Lambda=\Delta_A$. We define the {\it rank} 
of $\Lambda=\Delta_A$ simply to be the rank 
of the matrix $A$. 

\item For any $n\geq 1$, $\Lambda\in\bD_2[n]$ is 
said to be {\it full rank} if $\Lambda$ 
has rank $n$. The set of full rank elements 
of $\bD_2[n]$ will be denoted by 
$\bD_2^\circ [n]$. 

\item For any $r\geq 1$, we set 
\begin{align}
\Delta_r\!:=\sum_{i=1}^r D_i^2.
\end{align}
Note that $\Delta_r$ is a full rank element 
in $\bD_2[r]$ but not in $\bD_2[n]$ 
for any $n>r$.\\
\end{enumerate} 

For $U\in GL(n, \bC)$, we define 
\begin{align}
\Phi_U\!:\bar{\cA_n} & \to\quad \bar{\cA_n} \\
            P(z)&\to P(U^{-1}z) \nno  
\end{align}
and 
\begin{align}
\Psi_U\!:\cD[n]  \quad &  \to \quad\quad  \cD[n] \\
          \Lambda \quad & \to \Phi_U \circ \Lambda \circ  \Phi_U^{-1} \nno  
\end{align}

It is easy to see that, $\Phi_U$ (resp.\,$\Psi_U$)
is an algebra automorphism of $\cA_n$ 
(resp.\,$\cD[n]$). Moreover,  
the following standard facts 
are also easy to check directly.

\begin{lemma}\label{L2.1}
$(a)$ For any $U=(u_{ij}) \in GL(n, \bC)$, $P(z)\in \bar{\cA}_n$ 
and $\Lambda\in \bD[n]$, 
we have
\begin{align}\label{L2.1-e1}
\Phi_U (\Lambda P)= \Psi_U (\Lambda) \Phi_U(P).
\end{align}

$(b)$ For any $1\leq i\leq n$ and $f(z)\in \cA_n$ 
we have
\begin{align*}
\Psi_U (D_i) &=\sum_{j=1}^n u_{ji}D_j,\\
\Psi_U (f(D)) &=f(U^\tau D).
\end{align*}

In particular, for any $A \in SM(n, \bC)$, we have
\begin{align}\label{L2.1-e2}
\Psi_U (\Delta_A)
=\Delta_{UAU^\tau}.
\end{align}
\end{lemma}

The following lemma will play a crucial role 
in our later arguments. Actually the lemma can be 
stated in a stronger form (see \cite{Hua}, for example)
which we do not need here. 

\begin{lemma}\label{I-r}
For any $A\in SM(n, \bC)$ of rank $r>0$, 
there exists $U \in GL(n, \bC)$
such that
\begin{align}\label{I-r-e}
A=U \begin{pmatrix} 
 I_{r\times r} & 0 \\
0 & 0
\end{pmatrix}
U^\tau
\end{align} 
\end{lemma}

Combining Lemmas \ref{L2.1} and \ref{I-r}, it is easy to 
see we have the following corollary.

\begin{corol}\label{C2.3}
For any $n\geq 1$ and $\Lambda, \Xi \in \bD_2[n]$ of same rank, 
there exists $U\in GL(n, \bC)$ such that $\Psi_U(\Lambda)=\Xi$.
\end{corol}

\subsection{The Vanishing Conjecture for the 2nd Order Homogeneous
Differential Operators with Constant Coefficients}
\label{S2.2}

In this subsection, we show that Conjecture \ref{LVC}, 
hence also $\JC$, is actually equivalent to 
$\VC$ or $\HVC$ for all 2nd order 
homogeneous  differential operators 
$\Lambda\in \bD_2$ (see Theorem \ref{GVC-MainThm}).
We also show that the Laplace operators are not the 
only choices in the study of $\VC$ or $\JC$ 
(see Proposition \ref{MainPropo-3} and 
Example \ref{Mo-Chois}).  

First, let us point out that 
$\VC$ fails badly for differential operators
with non-constant coefficients. The following 
counter-example was given by M. de Bondt \cite{Bondt}.

\begin{exam}\label{Bondt} 
Let $x$ be a free variable and 
$\Lambda=x \frac{d^2}{d x^2}$. Let $P(x)=x$. Then 
one can check inductively that $P(x)$ is $\Lambda$-nilpotent, 
but $\Lambda^m P^{m+1}\neq 0$ for any $m\geq 1$.
\end{exam}

\begin{lemma}\label{L2.1-1}
For any $\Lambda\in \cD[n]$,  $U \in GL(n, \bC)$,  
$A \in SM(n, \bC)$ and $P(z)\in \bar{\cA}_n$, 
we have
\begin{enumerate}
\item[$(a)$]  $P(z)$ is $\Lambda$-nilpotent iff $\Phi_U (P)$ 
is $\Psi_U (\Lambda)$-nilpotent. In particular,
$P(z)$ is $\Delta_A$-nilpotent iff $\Phi_U (P)=P(U^{-1}z)$ 
is $\Delta_{UAU^\tau}$-nilpotent.

\item[$(b)$] $\VC[n]$ $(resp.\, \HVC[n])$ holds for $\Lambda$ 
iff it holds for $\Psi_U(\Lambda)$.
In particular, $\VC[n]$ $(resp.\,\HVC[n])$ 
holds for $\Delta_A$ iff 
it holds for $\Delta_{UAU^\tau}$.
\end{enumerate}
\end{lemma}

\pf Note first that, for any $m, k\geq 1$, 
we have
\begin{align*}
\Phi_U \left( \Lambda^m P^k \right)&=(\Phi_U \Lambda^m \Phi_U^{-1})\,  
\Phi_U P^k\\
&=(\Phi_U \Lambda \Phi_U^{-1})^m (\Phi_U P )^k \nno \\
&= [\Psi_U(\Lambda)]^m (\Phi_U P )^k. \nno
\end{align*}

When $\Lambda=\Delta_A$, by Eq.\,$(\ref{L2.1-e2})$, 
we further have
\begin{align*}
\Phi_U \left(\Delta_A^m P^k\right) = \Lambda_{UAU^\tau}^m (\Phi_U P)^k. \nno
\end{align*}

Since $\Phi_U$ (resp.\,$\Psi_U$) 
is an automorphism of $\bar{\cA}_n$ 
(resp.\,$\cD[n]$), 
it is easy to check directly that 
both $(a)$ and $(b)$ follow from the 
equations above.
\epfv

Combining the lemma above with Corollary \ref{C2.3}, 
we immediately have the following corollary.

\begin{corol}\label{C-nn}
Suppose $\HVC[n]$ $($resp.\,$\VC[n]$$)$ holds for 
a differential operator $\Lambda\in \bD_2 [n]$ 
of rank $r\geq 1$. 
Then $\HVC[n]$ $($resp.\,$\VC[n]$$)$ 
holds for all differential operators 
$\Xi\in \bD_2[n]$ of rank $r$.
\end{corol}

Actually we can derive much more (as follows) from 
the conditions in the corollary above.

\begin{propo}\label{MainPropo-1}
$(a)$  Suppose $\HVC[n]$ holds 
for a full rank $\Lambda \in \bD_2^\circ[n]$. 
Then, for any $k \leq n$, $\HVC[k]$ holds 
for all full rank $\Xi \in \bD_2^\circ[k]$.

$(b)$ Suppose $\VC[n]$ holds 
for a full rank $\Lambda \in \bD_2^\circ [n]$.
Then, for any $m \geq n$, $\VC[m]$ holds 
for all  $\Xi \in \bD_2[m]$ of rank $n$. 
\end{propo}

\pf Note first that, the cases $k=n$ in $(a)$ and 
$m=n$ in $(b)$ follow directly from Corollary \ref{C-nn}.
So we may assume $k<n$ in $(a)$ and $m>n$ in $(b)$.
Secondly, by Corollary \ref{C-nn}, it will be enough 
to show $\HVC[k]$ $(k<n)$ holds for $\Delta_k$ for $(a)$
and $\VC[m]$ $(m>n)$ holds for $\Delta_n$ for $(b)$.

$(a)$ Let $P\in \cA_k$ a homogeneous 
$\Delta_k$-nilpotent polynomial.
We view $\Delta_k$ and $P$ as elements 
of $\bD_2[n]$ and $\cA_n$, respectively.
Since $P$ does not 
depend on $z_i$ $(k+1\leq i\leq n)$, 
for any $m, \ell \geq 0$,
we have
\begin{align*}
\Delta_k^m P^\ell=\Delta_n^m P^\ell.
\end{align*}
Hence, $P$ is also $\Delta_n$-nilpotent. 
Since $\HVC[n]$ holds for $\Delta_n$ (as pointed out 
at the beginning of the proof), 
we have $\Delta_k^m P^{m+1} =\Delta_n^m P^{m+1}=0$ 
when $m>>0$. Therefore, $\HVC[k]$ holds for $\Delta_k$. 

$(b)$ Let $K$ be the rational function field 
$\bC(z_{n+1}, \dots, z_m)$.
We view $\cA_m$ as a subalgebra of 
the polynomial algebra $K[z_1, \dots, z_n]$ 
in the standard way. Note that the differential operator
$\Delta_n=\sum_{i=1}^n D_i^2$ of $\cA_m$ extends 
canonically to a differential operator of 
$K[z_1, \dots, z_n]$ with constant 
coefficients.  

Since $\VC[n]$ holds for $\Delta_n$ 
over the complex field (as pointed out 
at the beginning of the proof), 
by Lefschetz's principle, 
we know that $\VC[n]$ also holds for 
$\Delta_n$ over the field $K$. 
Therefore, for any $\Delta_n$-nilpotent 
$P(z) \in \cA_m$, by viewing $\Delta_n$ 
as an element of $\bD_2(K[z_1, \dots, z_n])$ and 
$P(z)$ an element of $K[z_1, \dots, z_n]$ 
(which is still $\Delta_n$-nilpotent in the new setting),
we have $\Delta_n^k P^{k+1}=0$ when $k>>0$. 
Hence $\VC[m]$ holds for $P(z)\in \cA_m$ 
and $\Delta_n\in \bD_2[m]$. 
\epfv

\begin{propo}\label{MainPropo-2}
Suppose $\HVC[n]$ holds for a differential operator 
$\Lambda \in \bD_2[n]$ with rank $r<n$. Then, 
for any $k \geq r$, $\VC[k]$ holds for 
all $\Xi\in \bD_2[k]$ of rank $r$.
\end{propo}

\pf First, by Corollary \ref{C-nn}, we know  
$\HVC[n]$ holds for $\Delta_r$. 
To show Proposition \ref{MainPropo-2}, 
by Proposition \ref{MainPropo-1}, $(b)$, 
it will be enough to show that  
$\VC[r]$ holds for $\Delta_r$.

Let $P \in \cA_r \subset \cA_n$ be 
a $\Delta_r$-nilpotent polynomial. 
If $P$ is homogeneous, there is nothing to prove
since, as pointed out above, $\HVC[n]$ holds for $\Delta_r$.  
Otherwise, we homogenize $P(z)$ 
to $\wtilde P \in \cA_{r+1} \subseteq \cA_n$.
Since $\Delta_r$ is a homogeneous 
differential operator, it is easy to see that, 
for any $m, k\geq 1$, $\Delta_r^m P^k=0$ iff
$\Delta_r^m \wtilde P^k=0$. Therefore,  
$\wtilde P \in \cA_n$ is also 
$\Delta_r$-nilpotent when we view $\Delta_r$ as 
a differential operator of $\cA_n$. 
Since $\HVC[n]$ holds for $\Delta_r$, 
we have that $\Delta_r^m \wtilde P^{m+1}=0$ 
when $m>>0$. Then, by the observation above again, 
we also have $\Delta_r^m P^{m+1}=0$ when $m>>0$. 
Therefore, $\VC[r]$ holds for $\Delta_r$.
\epfv

Now we are ready to prove our main result of this section.

\begin{theo}\label{GVC-MainThm}
The following statements are equivalent to each other.
\begin{enumerate}
\item $\JC$  holds.
\item $\HVC[n]$ $(n\geq 1)$ hold for the Laplace operator $\Delta_n$. 
\item $\VC[n]$ $(n\geq 1)$ hold for the Laplace operator $\Delta_n$. 
\item $\HVC[n]$ $(n\geq 1)$ hold for all $\Lambda\in \bD_2[n]$.
\item $\VC[n]$ $(n\geq 1)$ hold for all $\Lambda\in \bD_2[n]$.
\end{enumerate}
\end{theo}

\pf First, the equivalences of $(1)$, $(2)$ and $(3)$ 
have been established in Theorem $7.2$ in \cite{HNP}.
While $(4)\Rightarrow (2)$,  $(5)\Rightarrow (3)$ 
and $(5)\Rightarrow (4)$ are trivial. Therefore, 
it will be enough to show $(3)\Rightarrow (5)$.

To show $(3)\Rightarrow (5)$, we fix any $n\geq 1$. 
By Corollary \ref{C-nn}, it will be enough
to show $\VC[n]$ holds for $\Delta_r$ $(1\leq r\leq n)$.
But under the assumption of $(3)$ (with $n=r$), 
we know that $\VC[r]$ holds for $\Delta_r$. 
Then, by Proposition \ref{MainPropo-1}, $(b)$, 
we know $\VC[n]$ also holds for $\Delta_r$.
\epfv

Next, we show that, to study $\HVC$, equivalently 
$\VC$ or $\JC$, the Laplace operators are not 
the only choices, even though 
they are the best in many situations.

\begin{propo}\label{MainPropo-3}
Let $\{n_k\,| \, k\geq 1\}$ be a 
strictly increasing sequence 
of positive integers and 
$\{\Lambda_{n_k} \,|\, k\geq 1\}$ 
a sequence of differential 
operators in $\bD_2$ with 
$\mbox{rank\,}(\Lambda_{n_k})=n_k$ 
$(k\geq 1)$. Suppose that, for any $k\geq 1$,  
$\HVC[N_k]$ holds for $\Lambda_{n_k}$ 
for some $N_k\geq n_k$.  
Then, the equivalent 
statements in Theorem \ref{GVC-MainThm} 
hold. 
\end{propo}

\pf  We show, under the assumption in the proposition, 
the statement $(2)$ in Theorem \ref{GVC-MainThm} holds, i.e. 
for any $n\geq 1$, $\HVC[n]$ $(n\geq 1)$ holds 
for the Laplace operator $\Delta_n \in \bD_2[n]$. 

For any fixed $n\geq 1$, let $k\geq 1$ such that 
$n_k \geq n$. If $N_k=n_k$, then, by Proposition 
\ref{MainPropo-1}, $(a)$, we have  
$\HVC[n]$ $(n\geq 1)$ holds for the Laplace operator 
$\Delta_n \in \bD_2[n]$. 
If $N_k>n_k$, then, by Proposition \ref{MainPropo-2}, 
we know $\VC[n_k]$ (hence also $\HVC[n_k]$) 
holds for $\Delta_{n_k}$. 
Since $n_k\geq n$, by Proposition 
\ref{MainPropo-1}, $(a)$ again, 
we know $\HVC[n]$ does hold 
for the Laplace operator 
$\Delta_n$.
\epfv

\begin{exam}\label{Mo-Chois}
Besides the Laplace operators, by Proposition \ref{MainPropo-3},
the following sequences of differential operators are also 
among the most natural choices.

\begin{enumerate}
\item Let $n_k=k$ $(k\geq 2)$ $($or 
any other strictly increasing sequence of positive 
integers$)$. Let $\Lambda_{n_k}$ 
be the ``Laplace operator" 
with respect to the standard Minkowski metric 
of $\bR^{n_k}$. Namely, choose
\begin{align}\label{Minkowski}
\Lambda_{k}=D_1^2- \sum_{i=2}^{k} D_i^2.
\end{align}

\item Choose $n_k=2k$ $(k\geq 1)$ 
$($or any other strictly increasing sequence 
of positive even numbers$)$.
Let $\Lambda_{2k}$ be the ``Laplace operator" 
with respect to the standard symplectic metric 
on $\bR^{2k}$, i.e. choose 
\begin{align}\label{symplectic}
\Lambda_{2k}=\sum_{i=1}^k D_i D_{i+k}.
\end{align}

\item We may also choose the complex Laplace 
operators $\Delta_{\bar \p}$ instead of the real 
Laplace operator $\Delta$. More precisely, 
we choose $n_k=2k$ for any $k\geq 1$ and view 
the polynomial algebra of $w_i$ $(1\leq i\leq 2k)$ over 
$\bC$ as the polynomial algebra 
$\bC[z_i, \bar z_i\,|\, 1\leq i\leq k ]$ 
by setting $z_i=w_i+\sqrt{-1}\, w_{i+k}$ for any $1\leq i\leq k$.
Then, for any $k\geq 1$, we set \begin{align}
\label{cx-Delta}
\Lambda_{k}= \Delta_{\bar \p, k}\!:= 
\sum_{i=1}^k \frac{\p^2}{\p z_i\p \bar z_i}.
\end{align}

\item More generally, we may also choose 
$\Lambda_k=\Delta_{A_{n_k}}$, where 
$n_k\in \bN$ and $A_{n_k} \in SM(n_k, \bC)$  
$($not necessarily invertible$)$ $(k\ge 1)$
with strictly increasing ranks.
\end{enumerate}
\end{exam}

\renewcommand{\theequation}{\thesection.\arabic{equation}}
\renewcommand{\therema}{\thesection.\arabic{rema}}
\setcounter{equation}{0}
\setcounter{rema}{0}

\section{\bf Some Properties of $\Delta_A$-Nilpotent Polynomials}
\label{S3}

As pointed earlier in Section \ref{S1} 
(see page \pageref{rmk-HN}), 
for the Laplace operators $\Delta_n$ $(n\geq 1)$, 
the notion  $\Delta_n$-nilpotency coincides with the 
notion of Hessian nilpotency. HN (Hessian nilpotent) 
polynomials or formal power series, their associated symmetric 
polynomial maps and Conjecture \ref{LVC} 
have been studied in \cite{BE2}, \cite{BE3}, 
\cite{Burgers}--\cite{OP-HNP} and \cite{EZ}. 
In this section, we apply Corollary \ref{C2.3}, 
Lemma \ref{L2.1-1} and also Lefschetz's 
principle to transform some results obtained 
in the references above to certain results on
$\Lambda$-nilpotent $(\Lambda \in \bD_2)$ 
polynomials or formal power series, $\VC$ for $\Lambda$ 
and also associated polynomial maps. Another purpose 
of this section is to give a short survey on some results 
on HN polynomials and Conjecture \ref{LVC} 
in the more general setting of $\Lambda$-nilpotent polynomials 
and $\VC$ for differential operators 
$\Lambda \in \bD_2$.

In subsection  \ref{S3.1}, we 
transform some results in \cite{Burgers} and \cite{HNP} 
to the setting of $\Lambda$-nilpotent polynomials 
for $\Lambda\in \bD_2[n]$ of full rank 
(see Theorems \ref{PDE1}--\ref{Heat}). 
In subsection  \ref{S3.2}, we derive four criteria 
for $\Lambda$-nilpotency $(\Lambda\in \bD_2)$ 
(see Propositions \ref{A-Crit-0}, \ref{A-Crit-1}, 
\ref{A-Crit-3} and \ref{A-Crit-4}). The criteria 
in this subsection are transformations of the 
criteria of Hessian nilpotency derived in 
\cite{HNP} and \cite{OP-HNP}.

In subsection  \ref{S3.3}, we transform
some results in \cite{BCW}, \cite{Wa} and \cite{Y} 
on $\JC$; \cite{BE2} and \cite{BE3} on 
symmetric polynomial maps; \cite{HNP}, \cite{OP-HNP} 
and \cite{EZ} on HN polynomials to certain 
results on $\VC$ for $\Lambda\in \bD_2$. 
In subsection  \ref{S3.4}, 
we consider $\VC$ for high order differential 
operators with constant coefficients. 
We mainly focus on the differential operators 
$\Lambda=D^{\bf a}$ $({\bf a}\in \bN^n)$.
Surprisingly, $\VC$ for these operators is  
equivalent to a conjecture 
(see Conjecture \ref{LP-Conj}) 
on Laurent polynomials, which is similar 
to a non-trivial theorem (see Theorem \ref{DK1}) 
first conjectured by O. Mathieu  \cite{Mat} and 
later proved by J. Duistermaat 
and W. van der Kallen \cite{DK}.

\subsection{Associated Polynomial Maps and PDEs} \label{S3.1}

Once and for all in this section, we fix any $n\geq 1$ and 
$A\in SM(n, \bC)$ of rank $1\leq r\leq n$. 
We use $z$ and $D$, unlike we did before, 
to denote the $n$-tuples $(z_1, z_2, \dots, z_n)$ and 
$(D_1, D_2, \dots, D_n)$, respectively.
We define a $\bC$-bilinear form 
$\la\cdot, \cdot\ra_A$ by setting 
$\la u, v\ra_A\!:=u^\tau A v$ for 
any $u, v\in \bC^n$. Note that, when 
$A=I_{n\times n}$, the bilinear form 
defined above is just the standard 
$\bC$-bilinear form of $\bC^n$, which 
we also denote by $\la\cdot, \cdot\ra$.

By Lemma \ref{I-r}, we may write $A$ as 
in Eq.\,$(\ref{I-r-e})$. 
For any $P(z)\in \bar{\cA_n}$, we set 
\begin{align}\label{wdP}
\wtilde P(z)=\Phi_U^{-1} P(z)=P(Uz).
\end{align}

Note that, by Lemma \ref{L2.1}, $(b)$, we have
$\Psi_U^{-1}(\Delta_A)=\Delta_r$. 
By Lemma \ref{L2.1-1}, $(a)$,  
$P(z)$ is $\Delta_A$-nilpotent iff 
$\wtilde P(z)$ is $\Delta_r$-nilpotent. 

\begin{theo}\label{PDE1}
Let $t$ be a central parameter. For any $P(z)\in \cA_n$ 
with $o(P(z))\ge 2$ and $A\in SGL(n, \bC)$, set 
$F_{A, t}(z)\!:=z-tA \nabla P(z)$. Then 
\begin{enumerate}
\item[$(a)$] there exists a unique $Q_{A, t}(z)\in \bC[t][[z]]$ 
such that the formal inverse map $G_{A, t}(z)$ of $F_{A, t}(z)$
is given by 
\begin{align}\label{PDE1-e1}
G_{A, t}(z)=z+ tA\nabla Q_{A, t}(z).
\end{align}

\item[$(b)$] The $Q_{A, t}(z)\in \bC[t][[z]]$ in $(a)$ is the unique 
formal power series solution of the following Cauchy problem:
\begin{align}
\begin{cases}\label{PDE1-e2}
\frac{\p \, Q_{A, t}}{\p t}(z) =
\frac 12 \, \la \nabla Q_{A, t}, \nabla Q_{A, t} \ra_A, \\
Q_{A, t=0}(z)=P(z).
\end{cases}
\end{align}
\end{enumerate}
\end{theo}

\pf Let $\wtilde P$ as given in Eq.\,$(\ref{wdP})$ and set
\begin{align} \label{wdF}
\wtilde F_{A, t}(z)=z-t\nabla \wtilde P(z). 
\end{align}
By Theorem $3.6$ in \cite{Burgers}, we know the formal inverse map 
$\wtilde G_{A, t}(z)$ of $\wtilde F_{A, t}(z)$ is given by 
\begin{align}\label{PDE1-pe1}
\wtilde G_{A, t}(z)=z+t\nabla \wtilde Q_{A, t}(z),
\end{align}
where $\wtilde Q_{A, t}(z)\in \bC[t][[z]]$ is the 
unique formal power series solution of the 
following Cauchy problem:
\begin{align}
\begin{cases}\label{PDE1-pe2}
\frac{\p \, \wtilde Q_{A, t}}{\p t}(z)  =\frac 12\,
\la \nabla \wtilde Q_{A, t}, \nabla \wtilde Q_{A,t} \ra, \\
\wtilde Q_{A, t=0}(z)  =\wtilde P(z).
\end{cases}
\end{align}

From the fact that $\nabla \wtilde P(z)=(U^\tau \nabla P)(Uz)$, 
it is easy to check that
\begin{align}\label{PDE1-pe3}
(\Phi_U \circ \wtilde F_{A, t}\circ \Phi_U^{-1})(z)  
=z- t A \nabla P(z)=F_{A, t}(z), 
\end{align}
which is the formal inverse map of 
\begin{align}\label{PDE1-pe4}
(\Phi_U \circ \wtilde G_{A, t}\circ \Phi_U^{-1}) (z) =z+ t 
(U \nabla \wtilde Q_{A, t})(U^{-1}z).
\end{align}

Set 
\begin{align}\label{wdQ}
Q_{A, t}(z)\!:=\wtilde Q_{A, t}(U^{-1}z).
\end{align}
Then we have
\begin{align}
\nabla Q_{A, t}(z)&=(U^\tau)^{-1}(\nabla \wtilde Q_{A, t})(U^{-1}z),\nno \\
U^\tau \nabla Q_{A, t}(z)&=(\nabla \wtilde Q_{A, t})(U^{-1}z), \label{PDE1-pe5} \\
\intertext{Multiplying $U$ to the both sides of the equation above 
and noticing that $A=UU^\tau$ by Eq.\,$(\ref{I-r-e})$ since $A$ is of full rank, we get}
A\nabla Q_{A, t}(z)&=(U\nabla \wtilde Q_{A, t})(U^{-1}z).\label{PDE1-pe6}
\end{align}
Then, combining Eq.\,$(\ref{PDE1-pe4})$ and the equation above, 
we see the formal inverse 
$G_{A, t}(z)$ of $F_{A, t}(z)$
is given by 
\begin{align}
G_{A, t}(z)=(\Phi_U \circ \wtilde G_{A, t}\circ \Phi_U^{-1})(z)
=z+tA\nabla Q_{A, t}(z).
\end{align}
Applying $\Phi_U$ to Eq.\,$(\ref{PDE1-pe2})$ and 
by Eqs.\,$(\ref{wdQ})$, $(\ref{PDE1-pe5})$, we see that 
$Q_{A, t}(z)$ is the unique formal power series 
solution of the Cauchy problem Eq.\,$(\ref{PDE1-e2})$.
\epfv

By applying the linear automorphism $\Phi_U$ of $\bC[[z]]$ 
and employing a similar argument as in the proof of 
Theorem \ref{PDE1} above, we can generalize 
Theorems $3.1$ and $3.4$ in \cite{HNP}
to the following theorem on $\Delta_A$-nilpotent 
$(A\in SGL(n, \bC))$ formal power series.

\begin{theo}\label{PDE2}
Let $A$, $P(z)$ and $Q_{A, t}(z)$ as in Theorem \ref{PDE1}. 
We further assume $P(z)$ is $\Delta_A$-nilpotent.
Then, 
\begin{enumerate}
\item[$(a)$] $Q_{A, t}(z)$ is the unique 
formal power series solution of the following Cauchy problem:
\begin{align}
\begin{cases}
\frac{\p \, Q_{A, t}}{\p t}(z)=
\frac 14 \, \Delta_A Q_{A, t}^2,  \\
Q_{A, t=0}(z)=P(z).
\end{cases}
\end{align}
\end{enumerate}

\item[$(b)$] For any $k\geq 1$, we have 
\begin{align}
Q_{A, t}^k(z)=\sum_{m\geq 1} \frac{t^{m-1}}{2^m m!(m+k)!} \,
\Delta_A^m P^{m+1}(z).
\end{align}
\end{theo}

Applying the same strategy to Theorem $3.2$ in \cite{HNP}, 
we get the following theorem. 

\begin{theo}\label{Heat}
Let $A$, $P(z)$ and $Q_{A, t}(z)$ as in Theorem \ref{PDE2}.
For any non-zero $s\in \bC$, set
\BQn
V_{t, s}(z)\!:=\exp(sQ_t(z))=\sum_{k=0}^\infty \frac {s^k Q_t^k(z)}{k!}.
\EQn
Then, $V_{t, s}(z)$ is the unique formal power 
series solution of the following 
Cauchy problem of the heat-like equation:
\begin{align}\label{Cauchy-5}
\begin{cases}
\frac {\p V_{t,s}} {\p t}(z) = \frac 1{2s}\, \Delta_A V_{t, s}(z),\\
U_{t=0, s}(z)= \exp(s P(z)).
\end{cases}
\end{align}
\end{theo}

\subsection{Some Criteria of $\Delta_A$-Nilpotency}\label{S3.2}

In this subsection, with the notation and 
remarks fixed in the previous 
subsection in mind,
we apply the linear automorphism $\Phi_U$
to transform some criteria of Hessian 
nilpotency derived in \cite{HNP} and \cite{OP-HNP} 
to criteria of $\Delta_A$-nilpotency 
$( A\in SM(n, \bC) )$ 
(see Proposition \ref{A-Crit-0}, \ref{A-Crit-1}, 
\ref{A-Crit-3} and \ref{A-Crit-4} below).

\begin{propo}\label{A-Crit-0}
Let $A$ be given as in Eq.\,$(\ref{I-r-e})$. 
Then, for any $P(z)\in \cA_n$, 
it is $\Delta_A$-nilpotent iff  
the submatrix of $U^\tau (\Hes P)\, U$
consisting of the first $r$ rows and $r$ columns 
is nilpotent.

In particular, when $r=n$, 
i.e. $\Delta_A$ is full rank,  
any $P(z)\in \bD_2[n]$
is $\Delta_A$-nilpotent iff  
$U^\tau (\Hes P)\, U$ 
is nilpotent.
\end{propo}

\pf Let $\wtilde P(z)$ be as in Eq.\,$(\ref{wdP})$. Then, 
as pointed earlier, $P(z)$ is $\Delta_A$-nilpotent 
iff $\wtilde P(z)$ is $\Delta_r$-nilpotent. 

If $r=n$, then by Theorem \ref{Crit-1} ,
$\wtilde P(z)$ is $\Delta_r$-nilpotent iff
$\Hes \wtilde P(z)$ is nilpotent. But note that 
in general we have
\begin{align}\label{A-Crit-0-pe1}
\Hes \wtilde P(z)=\Hes P(Uz)=U^\tau [(\Hes P)(Uz)]\, U.
\end{align}
Therefore, $\Hes \wtilde P(z)$ is nilpotent iff 
$U^\tau [(\Hes P)(Uz)]\, U$ is nilpotent 
iff, with $z$ replaced by $U^{-1}z$, 
$U^\tau [(\Hes P)(z)]\, U$ is nilpotent.
Hence the proposition follows in this case.

Assume $r<n$. We view $\cA_r$ as a subalgebra of 
the polynomial algebra $K[z_1, \dots, z_r]$, where
$K$ is the rational field $\bC(z_{r+1},\dots, z_n)$.
By Theorem \ref{Crit-1}  and Lefschetz's principle, 
we know that $\wtilde P$ is $\Delta_r$-nilpotent iff
the matrix 
$
\left (\frac{\p^2\wtilde P}{\p z_i \p z_j} \right)_{1\leq i, j\leq r}
$
is nilpotent.

Note that the matrix 
$
\left (\frac{\p^2\wtilde P}{\p z_i \p z_j} \right)_{1\leq i, j\leq r}
$
is the submatrix of $\Hes \wtilde P(z)$ 
consisting of the first $r$ 
rows and $r$ columns. By Eq.\,$(\ref{A-Crit-0-pe1})$, 
it is also the submatrix of $U^\tau [\Hes P(Uz)]\,U$ 
consisting of the first $r$ 
rows and $r$ columns.
Replacing $z$ by $U^{-1}z$ 
in the submatrix above, we see  
$
\left (\frac{\p^2\wtilde P}{\p z_i \p z_j} \right)_{1\leq i, j\leq r}
$
is nilpotent iff the submatrix of 
$U^\tau[\Hes P(z)]\,U$ 
consisting of the first $r$ 
rows and $r$ columns is nilpotent. 
Hence the proposition follows.
\epfv

Note that, for any homogeneous quadratic polynomial
$P(z)=z^\tau B z$ with $B\in SM(n, \bC)$, 
we have $\Hes P(z)=2B$. Then, by Proposition \ref{A-Crit-0}, 
we immediately have the following corollary.

\begin{corol}\label{2-corol}
For any homogeneous  quadratic polynomial 
$P(z)=z^\tau B z$ with $B\in SM(n, \bC)$, 
it is $\Delta_A$-nilpotent iff the submatrix of \, 
$U^\tau B \, U$ consisting of the first 
$r$ rows and $r$ columns is nilpotent.
\end{corol}

\begin{propo}\label{A-Crit-1}
Let $A$ be given as in Eq.\,$(\ref{I-r-e})$. 
Then, for any $P(z)\in \bar{\cA}_n$ with $o(P(z))\geq 2$, 
 $P(z)$ is $\Delta_A$-nilpotent iff 
$\Delta_A^m P^m=0$ for any $1\leq m\leq r$.
\end{propo}

\pf Again, we let $\wtilde P(z)$ be 
as in Eq.\,$(\ref{wdP})$ and note that 
$P(z)$ is $\Delta_A$-nilpotent iff $\wtilde P(z)$ 
is $\Delta_r$-nilpotent. 

Since $r\leq n$. We view $\cA_r$ as a subalgebra of 
the polynomial algebra $K[z_1, \dots, z_r]$, where
$K$ is the rational field $\bC(z_{r+1},\dots, z_n)$.
By Theorem \ref{Crit-1} and Lefschetz's principle 
(if $r<n$), we have $\wtilde P(z)$ is $\Delta_r$-nilpotent iff
$\Delta_r^m \wtilde P^m=0$ for any $1\leq m\leq r$. 
On the other hand, by Eqs.\,$(\ref{L2.1-e1})$ and $(\ref{L2.1-e2})$, 
we have $\Phi_U \left(\Delta_r^m \wtilde P^m\right)
=\Delta_A^m P^m$ for any $m\ge 1$. Since $\Phi_U$ 
is an automorphism of $\cA_n$, we have that, 
$\Delta_r^m \wtilde P^m=0$ for any $1\leq m\leq r$
iff $\Delta_A^m P^m=0$ for any $1\leq m\leq r$.
Therefore, 
$\wtilde P(z)$ is $\Delta_A$-nilpotent iff 
$\Delta_A^m P^m=0$ for any $1\le m\le r$.
Hence the proposition follows. 
\epfv

\begin{propo}\label{A-Crit-3}
For any $A\in SGL(n, \bC)$ and any homogeneous 
$P(z)\in \cA_n$ of degree $d\geq 2$, we have, $P(z)$ is $\Delta_A$-nilpotent 
iff, for any $\beta\in \bC$, $(\beta_D)^{d-2} P(z)$ 
is $\Lambda$-nilpotent, where $\beta_D\!:=\la\beta, D\ra$.
\end{propo}

\pf Let $A$ be given as in Eq.\,$(\ref{I-r-e})$ and 
$\wtilde P(z)$ as in Eq.\,$(\ref{wdP})$.
Note that, $\Psi_U^{-1}(\Delta_A)=\Delta_n$ 
(for $\Delta_A$ is of full rank), and
$P(z)$ is $\Delta_A$-nilpotent 
iff $\wtilde P(z)$ is 
$\Delta_n$-nilpotent. 

Since $\wtilde P$ is also homogeneous  of degree $d\geq 2$, 
by Theorem \ref{Crit-1} in \cite{OP-HNP}, we know that, 
$\wtilde P(z)$ is $\Delta_n$-nilpotent iff, for any $\beta\in \bC^n$, 
$\beta_D^{d-2} \wtilde P$ is $\Delta_n$-nilpotent. 
Note that, from Lemma \ref{L2.1}, $(b)$, we have
\allowdisplaybreaks{
\begin{align*}
\Psi_U (\beta_D)&=\la \beta, U^\tau D\ra\\
&=\la U\beta, D\ra \\
&=(U\beta)_D,
\end{align*} }
and  
\begin{align}
\Phi_U (\beta_D^{d-2} \wtilde P)= \Psi_U (\beta_D)^{d-2} \Phi_U (\wtilde P)
=(U\beta)_D^{d-2} P.
\end{align}

Therefore, by Lemma \ref{L2.1-1}, $(a)$, 
 $\beta_D^{d-2} \wtilde P$ is $\Delta_n$-nilpotent 
iff $(U\beta)_D^{d-2} P$ is 
$\Delta_A$-nilpotent since $\Psi_U(\Delta_n)=\Delta_A$. 
Combining all equivalences above, we have $P(z)$
is $\Delta_n$-nilpotent iff, for any $\beta\in \bC^n$, 
$(U\beta)_D^{d-2} P$ is $\Delta_A$-nilpotent.
Since $U$ is invertible, when $\beta$ runs over $\bC^n$ 
so does $U\beta$. Therefore the proposition follows.
\epfv

Let $\{e_i\,|\, 1\leq i\leq n\}$ be the standard basis of $\bC^n$. Applying the proposition above to $\beta=e_i$ $(1\leq i\leq n)$, 
we have the following corollary.

\begin{corol}\label{C3.8}
For any homogeneous  $\Delta_A$-nilpotent polynomial $P(z)\in \cA_n$ 
of degree $d\ge 2$, $D_i^{d-2}P(z)$  
$(1\leq i\leq n)$ are also $\Delta_A$-nilpotent. 
\end{corol}

We think that Proposition \ref{A-Crit-3} and Corollary \ref{C3.8}
are interesting because, due to Corollary \ref{2-corol}, it is much 
easier to decide whether a quadratic form 
is $\Delta_A$-nilpotent or not.  

To state the next criterion, we need fix 
the following notation.

For any $A\in SGL(n, \bC)$, we let 
$X_A (\bC^n)$ be the set of isotropic vectors $u\in \bC^n$ 
with respect to the $\bC$-bilinear form $\la\cdot, \cdot\ra_A$. 
When $A=I_{n\times n}$, we also 
denote $X_A(\bC^n)$ simply by 
of $X(\bC^n)$. 

For any $\beta\in \bC^n$, we set 
$h_\alpha(z)\!:=\la\alpha, z\ra$. Then, 
by applying $\Phi_U$ to a well-known theorem on 
classical harmonic polynomials, which is the following 
theorem for $A=I_{n\times n}$ 
(see, for example, \cite{H} and \cite{T}), 
we have the following result on homogeneous 
$\Delta_A$-nilpotent polynomials.

\begin{theo}\label{T8.1.1}
Let $P$ be any homogeneous polynomial of degree 
$d\geq 2$ such that $\Delta_A P=0$. 
We have
\begin{align} \label{d-Form}
P(z)=\sum_{i=1}^k h_{\alpha_i}^d (z)
\end{align}
for some $k\ge 1$ and 
$\alpha_i\in X_A(\bC^n)$ $(1\leq i\leq k)$.
\end{theo}

Next, for any homogeneous  polynomial $P(z)$ of degree $d\ge 2$, we introduce
the following matrices:
\begin{align}
\Xi_P\!:=&\left(\la \alpha_i, \alpha_j \ra_A \right)_{k\times k}, \\
\Omega_P\!:=&\left(\la \alpha_i, \alpha_j \ra_A \, 
h_{\alpha_j}^{d-2}(z) \right)_{k\times k}.
\end{align}

Then, by applying $\Phi_U$ to Proposition $5.3$ in 
\cite{HNP} (the details will be omitted here), 
we have the following criterion 
of $\Delta_A$-nilpotency for homogeneous 
polynomials. 

\begin{propo}\label{A-Crit-4}
Let $P(z)$ be as given in Eq.\,$(\ref{d-Form})$.
Then $P(z)$ is $\Delta_A$-nilpotent iff 
the matrix $\Omega_P$ is nilpotent.
\end{propo}

One simple remark on the criterion above is as follows.

Let $B$ be the $k\times k$ diagonal matrix with 
$h_{\alpha_i} (z)$ $(1\leq i\leq k)$ as the $i^{th}$  
diagonal entry. For any $1\leq j\leq k$, set
\begin{align}
\Omega_{P; j}\!:= B^j \Xi_P B^{d-2-j}=
(h_{\alpha_i}^j \la \alpha_i, \alpha_j \ra h_{\alpha_j}^{d-2-j}).
\end{align}

Then, by repeatedly applying the fact that, 
for any $C, D\in M(k, \bC)$, $CD$ is nilpotent iff so is $DC$, 
it is easy to see that Proposition \ref{A-Crit-4} can also 
be re-stated as follows.

\begin{corol}
Let $P(z)$ be given by Eq.\,$(\ref{d-Form})$ with $d\geq 2$. 
Then, for any $1\leq j\leq d-2$ and $m\geq 1$, 
$P(z)$ is $\Delta_A$-nilpotent iff 
the matrix $\Omega_{P; j}$ is nilpotent.
\end{corol}

Note that, when $d$ is even, we may choose $j=(d-2)/2$. 
So $P$ is $\Delta_A$-nilpotent iff the symmetric matrix 
\begin{align}\label{(d-2)}
\Omega_{P; (d-2)/2}=( h_{\alpha_i}^{(d-2)/2} \la \alpha_i, \alpha_j \ra_A 
h_{\alpha_j}^{(d-2)/2})
\end{align}
is nilpotent.

\subsection{Some Results on the Vanishing Conjecture 
of the 2nd Order Homogeneous Differential Operators 
with Constants Coefficients}\label{S3.3}

In this subsection, we transform some known results of 
$\VC$ for the Laplace operators $\Delta_n$ $(n\ge 1)$ 
to certain results on $\VC$ for $\Delta_A$ 
$(A\in SGL(n, \bC))$.

First, by Wang's theorem \cite{Wa}, we know that 
$\JC$  holds for any polynomial maps $F(z)$ 
with $\deg F\leq 2$. Hence, $\JC$  also holds 
for symmetric polynomials $F(z)=z-\nabla P(z)$ 
with $P(z)\in \bC[z]$ 
of degree $d\le 3$. By the equivalence of $\JC$ and  
$\VC$ for the Laplace operators 
established in \cite{HNP}, we know $\VC$ holds 
if $\Lambda=\Delta_n$ and $P(z)$ is a HN 
polynomials of degree $d \le 3$. 
Then, applying the linear automorphism 
$\Phi_U$, we have the following 
proposition.

\begin{theo}
For any $A\in SGL(n, \bC)$ and 
$\Delta_A$-nilpotent $P(z)\in \cA_n$  
$($not necessarily homogeneous$)$ 
of degree $d\leq 3$, 
we have $\Lambda^m P^{m+1}=0$ when $m>>0$, 
i.e. $\VC[n]$ holds for $\Lambda$ and $P(z)$.
\end{theo} 

Applying the classical homogeneous  reduction on $\JC$  
(see \cite{BCW}, \cite{Y}) to associated symmetric maps, 
we know that, to show $\VC$ for $\Delta_n$ $(n\ge 1)$, 
it will be enough to consider only homogeneous HN polynomials 
of degree $4$. Therefore, by applying the linear automorphism 
$\Phi_U$ of $\cA_n$, we have the same reduction 
for $\HVC$ too.

\begin{theo}\label{degree-4}
To study $\HVC$ in general, it will 
be enough to consider only homogeneous 
$P(z)\in \cA$ of degree $4$.
\end{theo}

In \cite{BE2} and \cite{BE3} it has been shown 
that $\JC$ holds for symmetric maps $F(z)=z-\nabla P(z)$ 
$(P(z)\in \cA_n)$ if the number of variables $n$ 
is less or equal to $4$, or $n=5$ and $P(z)$ is 
homogeneous. By the equivalence of $\JC$ for symmetric polynomial 
maps and $\VC$ for the Laplace operators 
established in \cite{HNP}, and 
Proposition \ref{MainPropo-2} 
and Corollary \ref{C-nn}, 
we have the following 
results on $\VC$ and $\HVC$.

\begin{theo}\label{BE2-3}
$(a)$ For any $n\geq 1$, $\VC[n]$ holds for any
$\Lambda \in \bD_2$ of rank $1\leq r\leq 4$. 

$(b)$ $\HVC[5]$ holds for any 
$\Lambda \in \bD_2[5]$.
\end{theo} 

Note that the following vanishing properties 
of HN formal power series 
have been proved in Theorem $6.2$ in \cite{HNP} 
for the Laplace operators $\Delta_n$ $(n\geq 1)$. 
By applying the linear automorphism $\Phi_U$, 
one can show it also holds for 
any $\Lambda$-nilpotent $(\Lambda\in \bD_2)$ 
formal power series. 

\begin{theo}\label{k-1}
Let $\Lambda\in \bD_2[n]$ and 
$P(z)\in \bar{\cA}_n$ be $\Lambda$-nilpotent 
with $o(P)\geq 2$.  The following statements are equivalent.
\begin{enumerate}
\item[(1)] $\Lambda^m P^{m+1}=0$ when $m>>0$.
\item[(2)] There exists $k_0\geq 1$ such that 
$\Lambda^m P^{m+k_0}=0$ when $m>>0$.
\item[(3)] For any fixed $k\geq 1$, 
$\Lambda^m P^{m+k}=0$ when $m>>0$.
\end{enumerate}
\end{theo}

By applying the linear automorphism 
$\Phi_U$, one can transform Theorem $1.5$ 
in \cite{EZ} on $\VC$ of the Laplace operators
to the following result on $\VC$ 
of $\Lambda\in \bD_2$.
 
\begin{theo}\label{GGVC}
Let $\Lambda\in \bD_2[n]$ and 
$P(z)\in \bar{\cA}_n$ any $\Lambda$-nilpotent polynomial
with $o(P)\geq 2$.  Then $\VC$ holds for $\Lambda$ and $P(z)$ 
iff, for any $g(z)\in \cA_n$, we have $\Lambda^m (g(z)P^m)=0$ 
when $m>>0$.
\end{theo}

In \cite{EZ}, the following theorem has also 
been proved for $\Lambda=\Delta_n$. 
Next we show it is also true in general.

\begin{theo}\label{EZ}
Let $A\in SGL(n, \bC)$ 
and $P(z)\in \cA_n$ a homogeneous $\Delta_A$-nilpotent polynomial 
with $\deg P \ge 2$. Assume that 
$\sigma_{A^{-1}}(z)\!:=z^\tau A^{-1}z$ 
and the partial derivatives $\frac{\p P}{\p z_i}$ $(1\leq i\leq n)$ 
have no non-zero common zeros. Then $\HVC[n]$ holds 
for $\Delta_A$ and $P(z)$.

In particular, if the projective subvariety 
determined by the ideal $\la P(z)\ra $ of $\cA_n$
is regular, $\HVC[n]$ holds for $\Delta_A$ and $P(z)$.
\end{theo}

\pf Let $\wtilde P$ as given in 
Eq.\,$(\ref{wdP})$. By Theorem $1.2$ in \cite{EZ}, 
we know that, when $\sigma_2(z)\!:=\sum_{i=1}^n z_i^2$ 
and the partial derivatives $\frac{\p \wtilde P}{\p z_i}$ 
$(1\leq i\leq n)$ have no non-zero common zeros, 
$\HVC[n]$ holds for $\Delta_n$ and $\wtilde P$. 
Then, by Lemma \ref{L2.1-1}, $(b)$, $\HVC[n]$ 
also holds for $\Delta_A$ and $P$.

But, on the other hand, since $U$ is invertible and, 
for any $1\leq i\leq n$,
$$
\frac{\p \wtilde P}{\p z_i}=\sum_{j=1}^n u_{ji}
\frac{\p P}{\p z_j}(Uz), 
$$
$\sigma_2(z)$ and $\frac{\p \wtilde P}{\p z_i}$ $(1\leq i\leq n)$ 
have no non-zero common zeros iff 
$\sigma_2(z)$ and $\frac{\p  P}{\p z_i}(Uz)$ $(1\leq i\leq n)$ 
have no non-zero common zeros, and iff, with $z$ replaced 
by $U^{-1}z$, $\sigma_2(U^{-1}z)=\sigma_{A^{-1}}(z)$ 
and $\frac{\p P}{\p z_i}(z)$ $(1\leq i\leq n)$ 
have no non-zero common zeros. 
Therefore, the theorem holds.
\epfv

\subsection{The Vanishing Conjecture for Higher 
Order Differential Operators with Constant Coefficients}
\label{S3.4}

Even though the most interesting case of 
$\VC$ is for $\Lambda\in \bD_2$, at least when 
$\JC$  is concerned, the case of $\VC$ for 
higher order differential operators with constant 
coefficients is also interesting and non-trivial. 
In this subsection, we mainly discuss $\VC$ 
for the differential operators $D^{\bf a}$ $({\bf a}\in \bN^n)$. 
At the end of this subsection, we also recall a result 
proved in \cite{OP-HNP} which says that, 
when the base field has characteristic $p>0$, $\VC$, 
even under a weaker condition, 
actually holds for any differential operator $\Lambda$ 
(not necessarily with constant coefficients).  

Let $\beta_j \in \bC^n$ $(1 \leq j \leq \ell)$ be linearly 
independent and set $\delta_j\!:=\la \beta_j, D\ra$.
Let $\Lambda=\prod_{j=1}^\ell \delta_j^{a_i}$ 
with $a_j \ge 1$ $(1\leq j\leq \ell)$. 

When $\ell=1$, $\VC$ for $\Lambda$ can be 
proved easily as follows. 

\begin{propo}\label{bD-1k}
Let $\delta\in \bD_1[z]$ and $\Lambda=\delta^k$ 
for some $k\geq 1$. Then 
\begin{enumerate}
\item[$(a)$] A polynomial $P(z)$ is $\Lambda$-nilpotent 
if $($and only if$)$ $\Lambda P=0$.

\item[$(b)$] $\VC$ holds for $\Lambda$.
\end{enumerate}
\end{propo} 

\pf Applying a change of variables, 
if necessary, we may assume 
$\delta=D_1$ and $\Lambda=D_1^k$. 

Let $P(z)\in \bCz$ such that 
$\Lambda P(z)=D_1^k P(z)=0$. 
Let $d$ be the degree of 
$P(z)$ in $z_1$. From the equation above, we have 
$k>d$. Therefore, for any $m\geq 1$, 
we have $km>dm$ which implies 
$\Lambda^m P(z)^m=D_1^{km} P^m(z)=0$. 
Hence, we have $(a)$.

To show $(b)$, let $P(z)$ 
be a $\Lambda$-nilpotent 
polynomial. By the same notation and 
argument above, we have $k>d$. 
Choose a positive integer $N>\frac{d}{k-d}$. 
Then, for any $m\geq N$, we have
$m>\frac{d}{k-d}$, 
which is equivalent to 
$(m+1)d < km$. Hence we have 
$\Lambda^m P(z)^{m+1}=D_1^{km} P^{m+1}(z)=0$.
\epfv

In particular, when $k=1$ in the proposition above, we  have 
the following corollary.

\begin{corol}\label{bD1}
$\VC$ holds for any  differential operator 
$\Lambda\in \bD_1$.
\end{corol}

Next we consider the case $\ell \ge 2$. Note that, 
when $\ell=2$ and $a_1=a_2=1$. $\Lambda\in \bD_2$ 
and has rank $2$. Then, by Theorem \ref{BE2-3}, 
we know $\VC$ holds for $\Lambda$. 

Besides the case above, $\VC$ for 
$\Lambda=\prod_{j=1}^\ell \delta_j^{a_i}$ 
with $\ell \ge 2$ seems to be non-trivial 
at all. Actually, we will show below, 
it is equivalent to a conjecture 
(see Conjecture \ref{LP-Conj}) 
on Laurent polynomials.

First, by applying a change of variables, 
if necessary, we may (and will) assume 
$\Lambda=D^{\bf a}$ with ${\bf a}\in (\bN^+)^\ell$. 
Secondly, note that, for any ${\bf b}\in \bN^n$ and 
$h(z)\in \bC[z]$, $D^{\bf b} h(z)=0$ iff 
the holomorphic part of the Laurent polynomial
$z^{-\bf b}h(z)$ is zero. 

Now we fix a $P(z)\in \bC[z]$ and set 
$f(z)\!:=z^{-\bf a}P(z)$. With the 
observation above, it is easy to see that, $P(z)$ 
is $D^{\bf a}$-nilpotent iff the holomorphic parts of 
the Laurent polynomials $f^m(z)$ $(m\geq 1)$ are all 
zero; and $\VC$ holds for $\Lambda$ and $P(z)$ 
iff the holomorphic part of $P(z)f^m(z)$ 
is zero when $m>>0$.  Therefore, 
$\VC$ for $D^{\bf a}$ can be restated as follows:\\

\underline{\bf Re-Stated $\VC$ for $\Lambda=D^{\bf a}$}: 
{\it Let $P(z)\in \cA_n$ and $f(z)$ as above. 
Suppose that, for any $m\ge 1$, 
the holomorphic part of the Laurent polynomial $f^m(z)$ 
is zero, then the holomorphic 
part of $P(z)f^m(z)$ equals to zero when $m>>0$.}\\

Note that the re-stated $\VC$ above is very similar to 
the following non-trivial theorem which was first 
conjectured by O. Mathieu \cite{Mat} and later proved 
by J. Duistermaat and W. van der Kallen \cite{DK}.

\begin{theo}\label{DK1}
Let $f$ and $g$ be Laurent polynomials in $z$. Assume that, 
for any $m\geq 1$, the constant term of $f^m$ is zero. Then
the constant term $g f^m$ equals to zero when $m>>0$.
\end{theo}

Note that, Mathieu's conjecture \cite{Mat} 
is a conjecture on all real 
compact Lie groups $G$, which is also mainly motivated 
by $\JC$. The theorem above is the special case of 
Mathieu's conjecture when $G$ the $n$-dimensional 
real torus. For other compact real Lie groups, 
Mathieu's conjecture seems to be 
still wide open.

Motivated by Theorem \ref{DK1}, 
the above re-stated $\VC$ for $\Lambda=D^{\bf a}$ 
and also the result on $\VC$ in Theorem \ref{GGVC}, 
we would like to propose the following conjecture 
on Laurent polynomials.

\begin{conj}\label{LP-Conj}
Let $f$ and $g$ be Laurent polynomials in $z$. Assume that, 
for any $m\geq 1$, the holomorphic part of $f^m$ is zero. Then
the holomorphic part $g f^m$ equals to zero when $m>>0$.
\end{conj}

Note that, a positive answer to the conjecture above will 
imply $\VC$ for $\Lambda=D^{\bf a}$ 
$({\bf a}\in \bN^n)$ by simply choosing $g(z)$ to be $P(z)$.

Finally let us to point out that, it is well-known that $\JC$  
does not hold over fields of finite characteristic 
(see \cite{BCW}, for example),
but, by Proposition $5.3$ in \cite{OP-HNP}, 
the situation for VC over 
fields of finite characteristic 
is dramatically different 
even though it is equivalent to $\JC$ over 
the complex field $\bC$.

\begin{propo}\label{pVC}
Let $k$ be a field of $char.\,p>0$ and 
$\Lambda$ any differential operator of $k[z]$. 
Let $f\in k[[z]]$. Assume that, 
for any $1\leq m\leq p-1$, 
there exists $N_m>0$ such that 
$\Lambda^{N_m} f^m=0$.
Then,  $\Lambda^{m} f^{m+1}=0$ 
when $m>>0$.
\end{propo}

From the proposition above, we immediately 
have the following corollary.

\begin{corol}\label{pVC-C1}
Let $k$ be a field of $char.\,p>0$. Then 

$(a)$ $\VC$ holds for any differential operator $\Lambda$ of $k[z]$. 

$(b)$ If $\Lambda$ strictly decreases the degree of polynomials. 
Then, for any polynomial $f\in k[z]$ $($not necessarily 
$\Lambda$-nilpotent$)$, we have $\Lambda^{m} f^{m+1}=0$ 
when $m>>0$.
\end{corol}

\renewcommand{\theequation}{\thesection.\arabic{equation}}
\renewcommand{\therema}{\thesection.\arabic{rema}}
\setcounter{equation}{0}
\setcounter{rema}{0}

\section{\bf A Remark on $\Lambda$-Nilpotent Polynomials
and Classical Orthogonal Polynomials}\label{S4}

In this section, we first in subsection  
\ref{S4.1} consider the ``formal" 
connection between 
$\Lambda$-nilpotent polynomials 
or formal power series and classical 
orthogonal polynomials, which has been 
discussed in Section \ref{S1} 
(see page \pageref{nil-orp1}). 
We then in subsection  \ref{S4.2} 
transform the isotropic properties 
of homogeneous HN polynomials proved 
in \cite{HNP} to isotropic properties 
of homogeneous $\Delta_A$-nilpotent 
$(A \in SGL(n, \bC))$ polynomials 
(see Theorem \ref{A-Isotropic} 
and Corollary \ref{A-Isotropic-C}). 
Note that, as pointed in Section \ref{S1}, 
the isotropic results in 
subsection  \ref{S4.2} can be 
understood as some natural 
consequences of the connection
of $\Lambda$-nilpotent polynomials
and classical orthogonal 
polynomials discussed 
in subsection  \ref{S4.1}.

\subsection{Some Classical Orthogonal Polynomials}\label{S4.1}

First, let us recall the definition of 
classical orthogonal polynomials. Note that, 
to be consistent with the tradition 
for orthogonal polynomials, we will 
in this subsection use $x=(x_1, x_2, \dots, x_n)$ 
instead of $z=(z_1, z_2, \dots, z_n)$ to denote free 
commutative variables.

\begin{defi}\label{Def-OP}
Let $B$ be an open set of $\bR^n$ and $w(x)$ 
a real valued function defined over $B$ such that 
$w(x)\geq 0$ for any $x\in B$ 
and $0<\int_B w(x)dx <\infty$. 
A sequence of polynomials 
$\{f_{\bf m}(x)\,|\, {\bf m}\in \bN^n\}$ is 
said to be {\it orthogonal} over $B$ 
if \\
$(1)$ $\deg f_{\bf m}=|{\bf m}|$ for any ${\bf m}\in \bN^n$. \\
$(2)$ $\int_B f_{\bf m}(x) f_{\bf k}(x) w(x) \, dx=0$ for any 
${\bf m}\neq {\bf k}\in \bN^n$.
\end{defi}

%

The function $w(x)$ is called the {\it weight function}.
When the open set $B\subset \bR^n$ and $w(x)$ are 
clear in the context, we simply call the polynomials 
$f_{\bf m}(x)$ $({\bf m}\in \bN^n)$ in the definition 
above {\it orthogonal polynomials}.
If the orthogonal polynomials 
$f_{\bf m}(x)$ $({\bf m}\in \bN^n)$ 
also satisfy $\int_B f_{\bf m}^2(x) w(x)dx=1$ 
for any ${\bf m}\in \bN^n$, we call $f_{\bf m}(x)$ 
$({\bf m}\in \bN^n)$ {\it orthonormal polynomials}.
Note that, in the one dimensional case $w(x)$ 
determines orthogonal polynomials over $B$ 
up to  multiplicative constants, i.e. if 
$f_m(x)$ $(m\geq 0)$ are orthogonal polynomials 
as in Definition \ref{Def-OP}, then, 
for any $a_m\in \bR^{\times}$ $(m\geq 0)$, 
$a_m f_m$ $(m\geq 0)$ are also orthogonal over 
$B$ with respect to the weight function 
$w(x)$. 

The most natural way to construct orthogonal 
or orthonormal sequences is: first to 
list all monomials in an order such that 
the degrees of monomials are 
non-decreasing; and then to apply Gram-Schmidt 
procedure to orthogonalize or orthonormalize 
the sequence of monomials. But, 
surprisingly, most of classical 
orthogonal polynomials can also be 
obtained by the so-called Rodrigues' 
formulas.

We first consider orthogonal polynomials in 
one variable.\\ 

{\bf Rodrigues' formula}: {\it Let 
$f_m(x)$ $(m\geq 0)$ be the orthogonal polynomials 
as in Definition \ref{Def-OP}. Then, there exist 
a function $g(x)$ defined over $B$ 
and non-zero constants 
$c_m \in \bR$ $(m\geq 0)$ 
such that }
\begin{align}\label{Rodrigues}
f_m(x)=c_m w(x)^{-1} \frac {d^m}{dx^m} (w(x)g^m(x)). \\ \nno
\end{align}

Let $P(x)\!:=g(x)$ and $\Lambda\!:=w(x)^{-1} \left (\frac {d}{dx}\right) w(x)$,  
where, throughout this paper any polynomial or function 
appearing in a (differential) operator always means 
the multiplication operator by the polynomial 
or function itself. Then, 
by Rodrigues' formula above, we see that the 
orthogonal polynomials 
$\{f_m(x)\,|\, m\geq 0\}$ have the form
\begin{align}
f_m(x)=c_m \Lambda^m P^m(x),
\end{align}
for any $m\geq 0$.

In other words, all orthogonal polynomials 
in one variable, up to multiplicative constants, 
has the form $\{\Lambda^m P^m\,|\, m\ge 0\}$ 
for a single differential operator $\Lambda$ and 
a single function $P(x)$.

Next we consider some of the most well-known 
classical orthonormal polynomials 
in one variable. For more details 
on these orthogonal polynomials, 
see \cite{Sz}, \cite{AS}, \cite{DX}.

\begin{exam}\label{HP}
{\bf (Hermite Polynomials)} 

$(a)$ $B=\bR$ and the weight function $w(x)=e^{-x^2}$. 

$(b)$ Rodrigues' formula: 
\begin{align*}
H_m(x)=(-1)^m e^{x^2} ( \frac d{dx})^m e^{-x^2}. 
\end{align*}

$(c)$ Differential operator $\Lambda$ and polynomial $P(x)$:
\begin{align*}
\begin{cases}
\Lambda  = e^{x^2} ( \frac d{dx} ) e^{-x^2}=\frac d{dx}-2x, \\
P(x)=1,
\end{cases}
\end{align*}

$(d)$ Hermite polynomials in terms of $\Lambda$ and $P(x)$:
\begin{align*}
H_m(x)=(-1)^m \, \Lambda^m P^m(x).
\end{align*}
\end{exam}

\begin{exam}
{\bf (Laguerre Polynomials)} 

$(a)$ $B=\bR^+$ and 
$w(x)=x^\alpha e^{-x}$ $(\alpha>-1)$.

$(b)$ Rodrigues' formula: 
\begin{align*}
L_m^\alpha (x)=\frac 1{m!} x^{-\alpha} e^{x}(\frac d{dx})^m (x^{m+\alpha} e^{-x}).
\end{align*}

$(c)$ Differential operator $\Lambda$ and polynomial $P(x)$:
\begin{align*}
\begin{cases}
\Lambda_\alpha  = x^{-\alpha} e^{x}( \frac d{dx} ) (e^{-x} x^\alpha ) =\frac d{dx}+ (\alpha x^{-1}-1), \\
P(x)=x,
\end{cases}
\end{align*}

$(d)$ Laguerre polynomials in terms of $\Lambda$ and $P(x)$:
\begin{align*}
L_m(x)=\frac 1{m!}\,  \Lambda^m P^m(x).
\end{align*}
\end{exam}

\begin{exam}
{\bf (Jacobi Polynomials)} 

$(a)$ $B=(-1, 1)$ and 
$w(x)=(1-x)^\alpha (1+x)^\beta$, where $\alpha, \beta >-1$.

$(b)$ Rodrigues' formula: 
\begin{align*}
P_m^{\alpha, \beta} (x)=\frac {(-1)^m}{2^m m!}
(1-x)^{-\alpha} (1+x)^{-\beta} (\frac d{dx})^m (1-x)^{\alpha+m} (1+x)^{\beta+m} .
\end{align*}

$(c)$ Differential operator $\Lambda$ and polynomial $P(x)$:
\begin{align*}
\Lambda  &= (1-x)^{-\alpha} (1+x)^{-\beta}
(\frac d{dx}) (1-x)^{\alpha} (1+x)^{\beta} \\
&=\frac{d}{dx}-\alpha(1-x)^{-1}+\beta(1+x)^{-1},
\end{align*}
and
\begin{align*}
P(x)=1-x^2.
\end{align*}

$(d)$ Laguerre polynomials in terms of $\Lambda$ and $P(x)$:
\begin{align*}
P^{\alpha, \beta}_m (x)=\frac {(-1)^m}{2^m m!}\, \Lambda^m P^m(x).
\end{align*}
\end{exam}

A very important special family of Jacobi polynomials are the 
{\it Gegenbauer polynomials} which are obtained by setting 
$\alpha=\beta=\lambda-1/2$ for some $\lambda>-1/2$. 
Gegenbauer polynomials are also called 
{\it ultraspherical polynomials} in the literature.

\begin{exam}\label{GP}
{\bf (Gegenbauer Polynomials)} 

$(a)$ $B=(-1, 1)$ and 
$w(x)=(1-x^2)^{\lambda-1/2}$, where $\lambda >-1/2$.

$(b)$ Rodrigues' formula: 
\begin{align*}
P_m^\lambda (x)=\frac {(-1)^m}{2^m(\lambda+1/2)_m} 
(1-x^2)^{1/2-\lambda}(\frac d{dx})^m (1-x^2)^{m+\lambda-1/2}.
\end{align*}
where, for any $c\in \bR$ and $k\in \bN$, 
$(c)_k=c(c+1)\cdots (c+k-1)$.

$(c)$ Differential operator $\Lambda$ and polynomial $P(x)$:
\begin{align}\label{Ge-OP-e1}
\Lambda &= (1-x^2)^{1/2-\lambda} (\frac d{dx}) (1-x^2)^{\lambda-1/2}\\
&=\frac d{dx} - \frac {(2\lambda-1)\,x}{(1-x^2)}, \nno
\end{align}
and
\begin{align*}
P(x)= 1-x^2.
\end{align*}

$(d)$ Laguerre polynomials in terms of $\Lambda$ and $P(x)$:
\begin{align*}
P_m^{\lambda} (x)=\frac {(-1)^m}{2^m(\lambda+1/2)_m} \, \Lambda^m P^m(x).
\end{align*}
\end{exam}

Note that, for the special cases with $\lambda=0, 1, 1/2$, 
the Gegenbauer Polynomials
$P_m^\lambda(x)$ are called the 
{\it Chebyshev polynomial of the first kind, the second kind} 
and {\it Legendre polynomials}, 
respectively. Hence all these classical 
orthogonal polynomials  also have the form of $\Lambda^m P^m$ 
$(m\geq 0)$ up to some scalar multiple constants 
$c_m$ with $P(x)=1-x^2$ and the corresponding 
special forms of the differential operator 
$\Lambda$ in Eq.\,$(\ref{Ge-OP-e1})$.

\begin{rmk}
Actually, the Gegenbauer polynomials are more closely 
and directly related with $\VC$ in some 
different ways. See \cite{EC} for more discussions 
on connections of the Gegenbauer 
polynomials with $\VC$.
\end{rmk}

Next, we consider some classical orthogonal polynomials 
in several variables. We will see that they can also be 
obtained from certain sequences of the form 
$\{\Lambda^m P^m\,|\, m\geq 0\}$ in a 
slightly modified way. One remark is that,  
unlike the one-variable case, 
orthogonal polynomials in several variables
up to multiplicative constants
are not uniquely determined 
by weight functions.

The first family of classical orthogonal 
polynomials in several variables 
can be constructed by taking Cartesian 
products of orthogonal polynomials 
in one variable as follows.

Suppose $\{f_m\,|\, m\geq 0\}$ is a sequence of 
orthogonal polynomials in one variable, say as given 
in Definition \ref{Def-OP}. We fix any $n\geq 2$ and 
set  
\begin{align}
W(x)\!:=&\prod_{i=1}^n w(x_i),\label{Wx} \\
f_{\bf m} (x)\!:= &\prod_{i=1}^n f_{m_i}(x_i),\label{fm}
\end{align}
for any $x\in B^{\times n}$ and ${\bf m}\in \bN^n$.

Then it is easy to see that the sequence
$\{ f_{\bf m} (x) \,|\, {\bf m}\in \bN^n \}$
are orthogonal polynomials over $B^{\times n}$ 
with respect to the weight function
$W(x)$ defined above. 

Note that, by applying the construction above 
to the classical one-variable orthogonal polynomials 
discussed in the previous examples, one gets 
the classical {\it multiple Hermite Polynomials},
{\it multiple Laguerre polynomials},
{\it multiple Jacobi polynomials} and 
{\it multiple Gegenbauer polynomials}, 
respectively. 
 
To see that the multi-variable orthogonal polynomials constructed above 
can be obtained from a sequence of the form 
$\{\Lambda^m P^m(x)\,|\, m\ge 0\}$, we suppose 
$f_m$ $(m\geq 0)$ have Rodrigues' formula Eq.\,$(\ref{Rodrigues})$.
Let $s=(s_1, \dots, s_n)$ be $n$ central formal parameters and 
set 
\begin{align}
\Lambda_s \!:=& W(x)^{-1}\left(\sum_{i=1}^n 
s_i\frac{\p}{\p x_i}\right) W(x),\label{Ls} \\  
P(x)\!:=& \prod_{i=1}^n g (x_i). \label{Px}
\end{align}

Let $V_{\bf m}(x)$ $({\bf m}\in \bN^n)$ 
be the coefficient of $s^{\bf m}$ in 
$\Lambda_s^{|\bf m|} P^{|\bf m|}(x)$. 
Then, from Eqs.\,$(\ref{Rodrigues})$, 
$(\ref{Wx})$--$(\ref{Px})$, 
it is easy to check that, 
for any ${\bf m}\in \bN^n$, we have
\begin{align}
f_{\bf m}(x) =c_{\bf m}
\frac{ {\bf m}!}{|\bf m|!}
\, V_{\bf m}(x),
\end{align}
where $c_{\bf m}=\prod_{i=1}^n c_{m_i}$.

Therefore, we see that any multi-variable 
orthogonal polynomials constructed as above 
from Cartesian products of one-variable orthogonal polynomials 
can also be obtained from a single differential operator 
$\Lambda_s$ and a single function $P(x)$ 
via the sequence 
$\{\Lambda_s^m P^m \,|\, m\ge 0\}$.

\begin{rmk}
Note that, one can also take 
Cartesian products of different kinds of one-variable
orthogonal polynomials to create more 
orthogonal polynomials in several variables. 
By a similar argument as above, we see that all these 
multi-variable orthogonal polynomials can also 
be obtained similarly from a single sequence 
$\{\Lambda_s^m P^m \,|\, m\ge 0\}$.
\end{rmk}

Next, we consider the following two examples of 
classical multi-variable orthogonal polynomials 
which are not Cartesian products of one-variable
orthogonal polynomials. 

\begin{exam}\label{OPB}
{\bf (Classical Orthogonal Polynomials over Unit Balls)} 

$(a)$ Choose $B$ to be the open unit ball ${\mathbb B}^n$ 
of $\bR^n$ and the weight function
\begin{align*}
W_\mu (x)=(1-||x||^2)^{\mu - 1/2},
\end{align*}
where $||x||=\sum_{i=1}^n x_i^2$ 
and $\mu > 1/2$.

$(b)$ Rodrigues' formula: For any ${\bf m}\in \bN^n$, set
\begin{align*}
U_{\bf m}(x)\!:=\frac{ (-1)^{\bf m} (2\mu)_{|\bf m|} }
{2^{|\bf m|}{\bf m}! (\mu +1/2)_{|\bf m|} }
\frac{\p^{|\bf m|}}{ \p x_1^{m_1} \cdots \p x_n^{m_n}} 
(1-||x||^2)^{|\bf m|+\mu -1/2}.
\end{align*}
Then, by Proposition $2.2.5$ in \cite{DX}, 
$\{U_{\bf m}(x)\, | \, {\bf m}\in \bN^n \}$
are orthonormal over ${\mathbb B}^n$ with respect to 
the weight function $W_\mu(x)$. 

$(c)$ Differential operator $\Lambda_s$ and polynomial $P(x)$: 
Let $s=(s_1, \dots, s_n)$ be $n$ 
central formal parameters and 
set 
\begin{align*}
\Lambda_s \!:=& W_\mu(x)^{-1}
\left(\sum_{i=1}^n s_i\frac{\p}{\p x_i}\right) W_\mu (x),\\  
P(x)\!:=& 1-||x||^2.
\end{align*}

Let $V_{\bf m}(x)$ $({\bf m}\in \bN^n)$ 
be the coefficient of $s^{\bf m}$ in 
$\Lambda_s^{|\bf m|}P^{|\bf m|}(x)$. 
Then from the Rodrigues type formula above, 
we have, for any ${\bf m}\in \bN^n$,
\begin{align*}
U_{\bf m}(x) =
\frac{ (-1)^{|\bf m|} (2\mu)_{|\bf m|}}{ 2^{|{\bf m}|} {|\bf m|}! 
(\mu+1/2)_{|\bf m|} }
\, V_{\bf m}(x).
\end{align*}

Therefore, the classical orthonormal polynomials
 $\{U_{\bf m}(x)\, | \, {\bf m}\in \bN^n \}$
over ${\mathbb B}^n$ can be obtained from 
a single differential operator $\Lambda_s$ and 
$P(x)$ via the sequence 
$\{\Lambda_s^m P^m \,|\, m\ge 0\}$.
\end{exam}

\begin{exam}\label{OPSX}
{\bf (Classical Orthogonal Polynomials over Simplices)} 

$(a)$ Choose $B$ to be the simplex 
\begin{align*}
T^n=\{ x\in \bR^n\, | \, \sum_{i=1}^n x_i<1; \,\, x_1, ..., x_n > 0 \}
\end{align*}
in $\bR^n$ and the weight function
\begin{align}
W_\kappa (x)=x_1^{\kappa_1}\cdots x_n^{\kappa_n}(1-|x|_1)^{ \kappa_{n+1} - 1/2},
\end{align}
where $\kappa_i>-1/2$ $(1\leq i\leq n+1)$ and $|x|_1=\sum_{i=1}^n x_i$.

$(b)$ Rodrigues' formula: For any ${\bf m}\in \bN^n$, set
\begin{align*}
U_{\bf m}(x)\!:=W_\kappa(x)^{-1} 
\frac{\p^{|\bf m|}}{ \p x_1^{m_1} \cdots \p x_n^{m_n}} 
\left( W_\kappa(x) (1-|x|_1)^{|\bf m|}\right).
\end{align*}
Then, $\{U_{\bf m}(x)\, | \, {\bf m}\in \bN^n \}$
are orthonormal over $T^n$ with respect to 
the weight function $W_\kappa (x)$. 
See Section $2.3.3$ of \cite{DX} 
for a proof of this claim.

$(c)$ Differential operator $\Lambda$ and polynomial $P(x)$: 
Let $s=(s_1, \dots, s_n)$ be $n$ 
central formal parameters and 
set 
\begin{align*}
\Lambda_s \!:=& W_\kappa (x)^{-1}
\left(\sum_{i=1}^n s_i\frac{\p}{\p x_i}\right) W_\kappa (x),\\  
P(x)\!:=& 1-|x|_1.
\end{align*}

Let $V_{\bf m}(x)$ $({\bf m}\in \bN^n)$ 
be the coefficient of $s^{\bf m}$ in 
$\Lambda_s^{|\bf m|}P^{|\bf m|}(x)$. 
Then from the Rodrigues type formula in $(b)$, 
we have, for any ${\bf m}\in \bN^n$,
\begin{align*}
U_{\bf m}(x) =
\frac{ {\bf m}!}{|\bf m|!}
\, V_{\bf m}(x).
\end{align*}

Therefore, the classical orthonormal polynomials
 $\{U_{\bf m}(x)\, | \, {\bf m}\in \bN^n \}$
over $T^n$ can be obtained from a single 
differential operator 
$\Lambda_s$ and a function $P(x)$ via the sequence 
$\{\Lambda_s^m P^m \,|\, m\ge 0\}$.
\end{exam}

\subsection{The Isotropic Property of $\Delta_A$-Nilpotent Polynomials}
\label{S4.2}

As discussed in Section \ref{S1}, 
the ``formal" connection 
of $\Lambda$-nilpotent 
polynomials with classical orthogonal 
polynomials predicts that $\Lambda$-nilpotent 
polynomials should be 
isotropic with respect 
to a certain $\bC$-bilinear 
form of $\cA_n$. In this subsection, 
we show that, for differential operators 
$\Lambda=\Delta_A$ $(A\in SGL(n, \bC))$,
this is indeed the case for any 
homogeneous $\Lambda$-nilpotent 
polynomials (see Theorem \ref{A-Isotropic} and 
Corollaries \ref{A-Isotropic-C}, \ref{A-Isotropic-C2}). \\

We fix any $n\geq 1$ and let 
$z$ and $D$ denote the n-tuples $(z_1, \dots, z_n)$
and $(D_1, D_2,\dots, D_n)$, respectively. 
Let $A\in SGL(n, \bC)$ and define the $\bC$-bilinear map 
\begin{align}
\{\cdot, \cdot\}_A: \cA_n  \times  \cA_n & \to \quad \cA_n \\
(f, \, g)  \quad  &\to  f(AD) g(z) , \nno
\end{align}

Furthermore, we also define a $\bC$-bilinear form 
\begin{align}\label{A-Bi-form}
(\cdot, \cdot)_A: \cA_n \times \cA_n & \to \quad \bC \\
(f,\, g) \quad &  \to \{f, g\}_{{}_A} |_{z=0}, \nno
\end{align}

It is straightforward to check that the $\bC$-bilinear 
form defined above is symmetric and its restriction 
on the subspace of homogeneous  polynomials of any fixed 
degree is non-singular.
Note also that, for any homogeneous polynomials 
$f, g\in \cA_n$ of the same degree, 
we have $\{f,  g\}_A=(f, g)_A$.

The main result of this subsection is the following theorem.

\begin{theo}\label{A-Isotropic}
Let $A\in SGL(n, \bC)$ and 
$P(z)\in \cA_n$ a homogeneous 
$\Delta_A$-nilpotent polynomial of degree $d\geq 3$. 
Let ${\mathcal I} (P)$ be the ideal of $\cA_n$
generated by $\sigma_{A^{-1}}(z)\!:=z^\tau A^{-1}z$ 
and $\fr {\p P}{\p z_i}$ $(1\leq i\leq n)$.
Then, for any $f(z)\in {\mathcal I} (P)$ 
and $m\geq 0$, we have 
\begin{align}\label{A-Isotropic-e}
\{f, \Delta_A^m P^{m+1}\}_A=f(AD) \, \Delta_A^m P^{m+1}=0.
\end{align}
\end{theo}

Note that, by Theorem $6.3$ in \cite{HNP}, we know that 
the theorem does hold when $A=I_n$ and $\Delta_A=\Delta_n$. 

\pf Note first that, elements of $\cA_n$ satisfying 
Eq.\,$(\ref{A-Isotropic-e})$ 
do form an ideal. Therefore, it will be enough to 
show $\sigma_{A^{-1}}(z)$ and $\fr {\p P}{\p z_i}$ 
$(1\leq i\leq n)$ satisfy Eq.\,$(\ref{A-Isotropic-e})$.
But Eq.\,$(\ref{A-Isotropic-e})$ 
for $\sigma_{A^{-1}}(z)$ simply follows 
the facts that $\sigma_{A^{-1}}(Az)=z^\tau Az$ and 
$\sigma_{A^{-1}}(AD)=\Delta_A$.

Secondly, by Lemma \ref{I-r}, we can write $A=UU^\tau$ for 
some $U=(u_{ij}) \in GL(n, \bC)$. Then, by Eq.\,$(\ref{L2.1-e2})$,  
we have $\Psi_U (\Delta_n)=\Delta_A$ or $\Psi_U^{-1} 
(\Delta_A) = \Delta_n$. 
Let $\wtilde P(z) \!:= \Phi_U^{-1} (P)=P(U z)$. 
Then by Lemma \ref{L2.1-1}, $(a)$, 
$\wtilde P$ is a homogeneous  $\Delta_n$-nilpotent polynomial, 
and by Eq.\,$(\ref{L2.1-e1})$, 
we also have
\begin{align}\label{A-Isotropic-pe1}
\Phi_U^{-1} (\Delta_A^m P^{m+1})=\Delta_n^m \wtilde P^{m+1}.
\end{align}

By Theorem $6.3$ in \cite{HNP}, 
for any $1\leq i\leq n$ and $m\geq 0$, we have, 
\begin{align*}
\frac{\p \wtilde P}{\p z_i}(D) \left(\Delta_n^m \wtilde P^{m+1}\right)=0
\end{align*}
Since 
\begin{align*}
\frac{\p \wtilde P}{\p z_i}(z)= \sum_{k=1}^n u_{ki} 
\frac{\p P}{\p z_k}(U z),
\end{align*}
we further have, 
\begin{align*}
\sum_{k=1}^n u_{ki} \frac{\p P}{\p z_k}(UD)
\left(\Delta_n^m \wtilde P^{m+1} \right)=0.
\end{align*}

Since $U$ is invertible, for any $1\leq i\leq n$, 
we have
\begin{align}\label{A-Isotropic-pe2}
\frac{\p P}{\p z_i}(U D)
\left(\Delta_n^m \wtilde P^{m+1} \right)=0.
\end{align}

Combining the equation above with Eq.\,$(\ref{A-Isotropic-pe1})$, 
we get
\begin{align} 
\frac{\p P}{\p z_i}(UD) \Phi_U^{-1} \left( \Delta_A^m  P^{m+1} \right)=0.\nno \\
\Phi_U^{-1}  (\Phi_U \frac{\p P}{\p z_i}(UD) \Phi_U^{-1} )
\left(\Delta_A^m  P^{m+1}\right)=0. \nno \\
 (\Phi_U \frac{\p P}{\p z_i}(UD) \Phi_U^{-1} )
\left(\Delta_A^m  P^{m+1}\right)=0. \label{A-Isotropic-pe3}
\end{align}

By Lemma \ref{L2.1}, $(b)$,  
Eq.\,$(\ref{A-Isotropic-pe3})$ and 
the fact that $A=UU^\tau$, we get
\begin{align*}
 \frac{\p P}{\p z_i}(U U^\tau D)
\left(\Delta_A^m  P^{m+1}\right)
=\frac{\p P}{\p z_i}(A D)
\left(\Delta_A^m  P^{m+1}\right)=0,
\end{align*}
which is Eq.\,$(\ref{A-Isotropic-e})$ for 
$\frac{\p P}{\p z_i}$ $(1\leq i\leq n)$.
\epfv

\begin{corol}\label{A-Isotropic-C}
Let $A$ be as in Theorem \ref{A-Isotropic} and 
$P(z)$ be a homogeneous $\Delta_A$-nilpotent polynomial 
of degree $d \geq 3$. Then, for any $m\geq 1$, 
$\Delta_A^m P^{m+1}$ is isotropic with respect to the 
$\bC$-bilinear form $(\cdot, \cdot)_A$, i.e.
\begin{align}
(\Delta_A^m P^{m+1}, \Delta_A^m P^{m+1})_A=0. 
\end{align}
In particular, we have $(P, P)_A=0$.
\end{corol}

\pf By the definition Eq.\,$(\ref{A-Bi-form})$ of 
the $\bC$-bilinear form $(\cdot, \cdot )_A$ and 
Theorem \ref{A-Isotropic}, it will be enough to show that 
$P$ and $\Delta_A^m P^{m+1}$ $(m\geq 1)$ belong to the ideal 
generated by the polynomials $\frac{\p P}{\p z_i}$ 
$(1\le i\le n)$ (here we do not need to consider 
the polynomial $\sigma_{A^{-1}}(z)$). 
But this statement has been proved 
in the proof of Corollary $6.7$ 
in \cite{HNP}. So we refer the reader to 
\cite{HNP} for a proof of the statement 
above. 
\epfv

Theorem \ref{A-Isotropic} and Corollary \ref{A-Isotropic-C} do not hold 
for homogeneous HN polynomials $P(z)$ of degree $d=2$. But, by applying 
similar arguments as in the proof of Theorem \ref{A-Isotropic} above
to Proposition $6.8$ in \cite{HNP}, one can show that 
the following proposition holds.

\begin{propo}\label{A-Isotropic-C2}
Let $A$ be as in Theorem \ref{A-Isotropic} and 
$P(z)$ a homogeneous $\Delta_A$-nilpotent polynomial 
of degree $d=2$. Let ${\mathcal J} (P)$ the ideal of $\bC[z]$
generated by $P(z)$ and $\sigma_{A^{-1}}(z)$.
Then, for any $f(z)\in {\mathcal J} (P)$ and $m\geq 0$, we have 
\begin{align}
\{f, \Delta_A^m P^{m+1}\}_A=f(AD) \, \Delta_A^m P^{m+1}=0.
\end{align}
In particular, we still have $(P, P)_A=0$.
\end{propo}

{\small \sc Department of Mathematics, Illinois State University,
Normal, IL 61790-4520.}

{\em E-mail}: wzhao@ilstu.edu.

\end{document}